\documentclass[12pt]{amsart}

\usepackage{amsmath}
\usepackage{amsfonts}
\usepackage{amssymb}
\usepackage{latexsym}
\usepackage{euscript}

\newtheorem{theorem}{Theorem}[section]
\newtheorem{proposition}[theorem]{Proposition}
\newtheorem{corollary}[theorem]{Corollary}

\newtheorem{definition}[theorem]{Definition}
\newtheorem{remark}[theorem]{Remark}
\theoremstyle{definition}
\newtheorem{example}[theorem]{Example}
\numberwithin{equation}{section}

   \def\sH{{\mathfrak H}}   
      
\def\sM{{\mathfrak M}}

      \def\dC{{\mathbb C}}

   \def\dN{{\mathbb N}}

   \def\dZ{{\mathbb Z}}

\def\cA{{\mathcal A}}      \def\cC{{\mathcal C}}
\def\cD{{\mathcal D}}      \def\cF{{\mathcal F}}
   \def\cH{{\mathcal H}}   
   \def\cK{{\mathcal K}}   \def\cL{{\mathcal L}}

\def\RE{{\rm Re\,}}
\def\IM{{\rm Im\,}}
\def\Ker{{\rm Ker\,}}
\def\wt{\widetilde}
\def\wh{\widehat}

\def\uphar{{\upharpoonright\,}}
\def\ovl{\overline}

\def\span{{\rm span\,}}
\def\ran{{\rm ran\,}}
\def\dom{{\rm dom\,}}

\begin{document}
\title[Non-self-adjoint Jacobi matrices]
{Non-self-adjoint Jacobi matrices with rank one imaginary part}
\author{
Yury~Arlinski\u{i}}
\address{Department of Mathematical Analysis \\
East Ukrainian National University \\
Kvartal Molodyozhny 20-A \\
Lugansk 91034 \\
Ukraine} \email{yma@snu.edu.ua}
\author{
Eduard~Tsekanovski\u{i}}
\address{Department of Mathematics, P.O. Box 2044 \\
Niagara University, NY 14109, USA } \email{tsekanov@niagara.edu}

 \maketitle
\begin{abstract} We develop direct and inverse spectral analysis
for finite and semi-infinite non-self-adjoint Jacobi matrices with a
rank one imaginary part. It is shown that given a set of $n$ not
necessarily distinct non-real numbers in the open upper (lower)
half-plane uniquely determines a $n\times n$ Jacobi matrix with a
rank one imaginary part having those numbers as its eigenvalues
counting
 multiplicity. An algorithm for reconstruction for such finite
Jacobi matrices is presented. A new model complementing the well
known Livsic triangular model for bounded linear operators with rank
one imaginary part is obtained. It turns out that the model operator
is a non-self-adjoint Jacobi matrix and it follows from the fact
that any bounded, prime, non-self-adjoint linear operator with rank
one imaginary part acting on some finite-dimensional (resp.,
separable infinite-dimensional Hilbert space) is unitary equivalent
to a finite (resp., semi-infinite) non-self-adjoint Jacobi matrix.
This result strengthens the classical Stone theorem established for
self-adjoint operators with simple spectrum. We establish the
non-self-adjoint analogs of the Hochstadt and Gesztesy--Simon
uniqueness theorems for finite Jacobi matrices with non-real
eigenvalues as well as an extension and refinement of these theorems
for finite non-self-adjoint tri-diagonal matrices to the case of
mixed eigenvalues, real and non-real. A unique Jacobi matrix,
unitarily equivalent to the operator of integration $(\cF
f)(x)=2\,i\,\int_{x}^l f(t)dt$ in the Hilbert space $\cL_2[0,l]$ is
found as well as spectral properties of its perturbations and
connections with well known Bernoulli numbers. We also give the
analytic characterization of the Weyl functions of dissipative
Jacobi matrices with a rank one imaginary part.

\end{abstract}

\tableofcontents

\section{Introduction}

Self-adjoint (or real) finite and semi-infinite Jacobi matrices of
the form
 \begin{equation}
\label{1} J=\begin{pmatrix} b_1 & a_1 & 0    & 0 & \cdot &
\cdot & \cdot \\
a_1 & b_2 & a_2 & 0 & \cdot &
\cdot & \cdot \\
0    & a_2 & b_3 & a_3 & \cdot &
\cdot & \cdot \\
\cdot & \cdot & \cdot & \cdot &
\cdot & \cdot & \cdot \\
\cdot & \cdot & \cdot & \cdot &
\cdot & \cdot & a_{n-1} \\
\cdot & \cdot & \cdot & \cdot & 0 & a_{n-1} & b_n
\end{pmatrix}
\end{equation}
and
\begin{equation}
\label{01}
  J=\begin{pmatrix} b_1 & a_1 & 0 &0   & 0 &
\cdot &
\cdot  \\
a_1 & b_2 & a_2 & 0 &0& \cdot &
\cdot  \\
0    & a_2 & b_3 & a_3 &0& \cdot &
\cdot   \\
\cdot & \cdot & \cdot & \cdot & \cdot & \cdot & \cdot
\end{pmatrix},
\end{equation}
where $a_k>0$, and $b_k$ are real numbers for all $k=1,2,\ldots$
play an important role in various problems of mathematical analysis
and theoretical and mathematical physics.  They appear as the
discrete analog of Sturm--Liouville operators, in inverse spectral
theory, in the study of the classical moment problem, and in the
investigation of completely integrable nonlinear lattices
\cite{Akh}, \cite{Ber}, \cite{Ber1}, \cite{Ch}, \cite{GS}, \cite{G},
  \cite{S}, \cite{Stone}, \cite{Tes}. More general tri-diagonal
matrices with complex entries (or complex Jacobi matrices) also
attracted much attention as a useful tool in the study of orthogonal
polynomials, in the theory of continued fractions, and in numerical
analysis \cite{Bec1}, \cite{Bec}, \cite{Wall}.

Let the linear space $\dC^n$ of columns be equipped by the usual
inner product
\[
(x,y)=\sum\limits_{k=1}^nx_k\overline{y_k}
\]
and let $l_2(\dN)$ be the Hilbert space of square summable
complex-valued sequences
$$x=\{x_1,x_2,\ldots,x_k,\ldots\}$$
which we consider as semi-infinite vector-columns with the inner
product given by
\[
(x,y)=\sum\limits_{k=1}^\infty x_k\overline{y_k}.
\]
Let
$$\delta_k:=(\underbrace {0,\ldots, 0,1}_k,0,\ldots)^T,\; k=1,2,\ldots.$$
 Then the vectors $\{\delta_k\}$ form an orthonormal basis in
$\dC^n$ (resp.,$l_2(\dN)$). An $n\times n$ complex Jacobi matrix $J$
determines a linear operator in the Hilbert space $\dC^n$ by means
of the matrix product $J\cdot x$. For the semi-infinite case we will
suppose in addition that
\begin{equation}
\label{Bound}
 \sup\limits_{k}\{|a_k|+|b_k|\}<\infty, k\in{\dN}.
\end{equation}
This  condition is necessary and sufficient \cite{AkhGl} for
boundedness of the Jacobi operator in $l_2(\dN)$ defined as $J\cdot
x$, where $J$ is a semi-infinite complex Jacobi matrix. Moreover,
under the condition $a_k\ne 0$ for all $k$, the vector $\delta_1$ is
cyclic for the Jacobi operator $J$. The classical Stone theorem
\cite{Akh}, \cite{AkhGl}, \cite{Stone} states that every
self-adjoint operator with simple spectrum in a separable Hilbert
space $\cH$ is unitarily equivalent to the operator determined by a
self-adjoint Jacobi matrix.

 In this paper we consider finite and semi-infinite tri-diagonal matrices of the form \eqref{1} and
\eqref{01} with
\begin{equation}
\label{333} \IM b_1>0,\;b_k={\bar b}_k\; (k=2, 3, ...),\;
a_k>0\;(k\in{\dN}).
\end{equation}
Such matrices determine bounded Jacobi operators $J$ in $\dC^n$ or
in $l_2(\dN)$ that posses the properties
 \[
\begin{split}
& \IM (Jx,x)\ge 0 \quad\mbox{for all}\quad x\in{\dC^n}~ or~ x\in l_2({\dN}),\\
&\ran(J-J^*)=\{\lambda \delta_1,\; \lambda\in\dC\}.
\end{split}
\]
In the following we will call such tri-diagonal matrices \textit{the
dissipative Jacobi matrices with rank one imaginary part}. Because
the vector $\delta_1$ is cyclic for $J$, the dissipative operator
$J$ is prime \cite{Br}, \cite{BrL}, \cite{L}, and therefore it has
no real eigenvalues.

M. S.~Livsic \cite{L} (see also \cite{BrL}) constructed a triangular
model for a bounded prime dissipative operator $A$ with a rank one
imaginary part. The method of construction is based on the
factorization of the characteristic function $W(z)$ of the operator
$A$. The model operator is in general a coupling of two triangular
operators. The first one is given by a finite or semi-infinite
triangular matrix of the form
\[
\dot{A}= \begin{pmatrix} \alpha_1+\frac{i}{2}\beta^2_1 &
i\beta_1\beta_2 & \cdot &
\cdot & \cdot & i\beta_1\beta_n \\
0 & \alpha_2+\frac{i}{2}\beta^2_2 &  \cdot &
\cdot & \cdot & i\beta_2\beta_n \\
\cdot & \cdot & \cdot & \cdot &
\cdot & \cdot  \\
0 & 0 & \cdot & \cdot & \cdot & \alpha_n+\frac{i}{2}\beta^2_n
\end{pmatrix}
\]
or
\[
\dot{A}=
\begin{pmatrix} \alpha_1+\frac{i}{2}\beta^2_1 & i\beta_1\beta_2 &
\cdot &
\cdot & \cdot & i\beta_1\beta_k&\cdot&\cdot&\cdot\\
0 & \alpha_2+\frac{i}{2}\beta^2_2 &  \cdot & \cdot & \cdot &
i\beta_2\beta_k&\cdot&\cdot&\cdot \\ \cdot & \cdot & \cdot & \cdot &
\cdot & \cdot& \cdot & \cdot & \cdot  \\
0 & 0 & \cdot & \cdot & \cdot & \alpha_k+\frac{i}{2}\beta^2_k& \cdot & \cdot & \cdot\\
\cdot&\cdot&\cdot&\cdot& \cdot & \cdot & \cdot& \cdot & \cdot
\end{pmatrix},
\]
where $\{\alpha_k\}_{k=1}^n$ are real numbers and
$\{\beta_k\}_{k=1}^n$, are positive numbers and
$z_k=\alpha_k+i(\beta^2_k/2), k=1, 2, ...$ are the non-real
eigenvalues of $\dot{A}$. The second one is the operator in the
Hilbert space $\cL_2[0,l]$ given by
\[
(\dot{B} f)(x)=\alpha(x)f(x)+i\int\limits_{x}^l f(t) dt,
\]
where $\alpha$ is a bounded nondecreasing right continuous function
on  $[0,l]$. The operator $\dot{B}$ has pure real spectrum which
coincides with the closure of the set $\{\alpha(x),\; x\in [0,l]\}$.

 One of the main results of the  paper states that \textit{a bounded
prime dissipative operator with a rank one imaginary part acting on
a separable Hilbert space is unitarily equivalent to the operator
determined by a dissipative Jacobi matrix with rank one imaginary
part. This result is a non-self-adjoint analog of the classical
Stone theorem, established for self-adjoint operators with simple
spectrum in \cite{Akh}, \cite{AkhGl}, \cite{Stone}. Thus,
dissipative Jacobi matrices with rank one imaginary parts provide
new models for the prime bounded linear operators with rank one
imaginary parts}. The entries of the corresponding Jacobi matrix can
be found using the continued fraction (J-fraction) expansion
\cite{Wall}
\[
M(z)= \frac{-1}{z-b_1}\;\raisebox{-3mm}{{\rm
+}}\;\frac{-a^2_1}{z-b_2}\;\raisebox{-3mm}{{\rm
+}}\;\frac{-a^2_2}{z-b_3} \;\raisebox{-3mm}{{\rm +}\,\ldots}\;
\raisebox{-3mm}{{\rm +}}\;\frac{-a^2_{n-1}}
{z-b_n}\;\raisebox{-3mm}{{\rm +}\,\dots},
\]
where  $W(z)$ is the characteristic function,
$M(z)=\frac{i}{\beta}(W(z)-1)$,
 and
 $\beta=\lim\limits_{z\to \infty}\left(iz(1-W(z))\right).$

By the Livsic theorem \cite{L}, a prime dissipative Volterra
operator $A$ acting on the Hilbert space $\cH$ with rank one
imaginary part such that
$$(A-A^*)h=2i\,l(h,e)e$$
 for all $h\in \cH$, $||e||=1$, is
unitarily equivalent to the integration operator of the form
\[
\c(F f)(x)=2i\int\limits_x^l f(t)dt
\]
acting on the Hilbert space $\cL_2[0,l]$. \textit {In this paper a
unique Jacobi matrix unitarily equivalent to the integration
operator $\cF$
 is found, and some connections with the well known Bernoulli numbers are established}.

 \textit {We also show that for $n$ not necessarily distinct given complex
numbers in the open upper half-plane, there exists a unique $n\times
n$ Jacobi matrix having those numbers as its eigenvalues counting
multiplicity.  Two algorithms for the reconstruction of a finite
 Jacobi matrix of the form \eqref{1} from its
 eigenvalues are presented}.
\textit {We obtain the following non-self-adjoint analogs of the 
 Hochstadt and Gesztesy--Simon uniqueness theorems for
self-adjoint finite Jacobi matrices: Let $J$ be an $n\times n$
Jacobi matrix satisfying conditions \eqref{333}. Suppose that the
eigenvalues $z_1,\ldots, z_k$ (taken from the open upper half-plane)
with their multiplicities $l_1, \ldots, l_k$ are known together with
$b_1,a_1,b_2,\dots, a_{n-r-1}, b_{n-r}$, where $r=l_1+l_2+\cdots+
l_k.$ Then $a_{n-r}, b_{n-r+1},$ $\ldots,a_{n-1}, b_n$ are uniquely
determined.} This as well as some \textit{extensions and refinements
of the Hochstadt and Gesztesy--Simon uniqueness theorems concerning
the case of dissipative but non prime tri-diagonal matrices with
rank one imaginary part are established.}


Our main tools are the Livsic characteristic function, its
linear-fractional transformation \cite{Br}, \cite{BrL}, \cite{L}
\cite{LYan}, and the Weyl function of a Jacobi matrix. Moreover, in
order to establish the above mentioned non-self-adjoint analogs of
the Hochstadt \cite{H1} and Gesztesy--Simon \cite{GS} uniqueness
theorems in the theory of inverse spectral problems we develop
further the approach Gesztesy--Simon used in \cite{GS}.

Our paper is organized as follows: In Section \ref{sec2} we present
some properties of prime, bounded, dissipative operators with a rank
one imaginary part and the Livsic theory of their characteristic
functions and triangular models. We define the Weyl function of such
operators and give its analytic characterizations. We emphasize that
the important role of the Weyl function (Weyl--Titchmarsh function)
in the spectral theory of self-adjoint differential and difference
operators and its applications to nonlinear equations is well-known
(see \cite{Akh},  \cite{Atk}, \cite{Ber}, \cite{Ber1}, \cite{GS},
\cite{GS1}, \cite{Tes}). In Section \ref{sec3} some basic properties
of complex Jacobi matrices and the corresponding Jacobi operators as
well as a survey of the inverse spectral problems for finite and
semi-infinite self-adjoint Jacobi matrices are given. We begin our
study of dissipative Jacobi matrices with a rank one imaginary part
in Section \ref{sec4}. In particular, the connection between $m_+$-
and $m_-$-functions of a finite dissipative Jacobi matrix is
established. In Section \ref{sec5} we present solutions of the
inverse spectral problems for finite dissipative Jacobi matrices.
The non-self-adjoint analogs of the Hochstadt and Gesztesy--Simon
theorems are presented in Subsections \ref{subsec5} and
\ref{subsec55}. The non-self-adjoint analog of Stone's theorem is
considered in Section \ref{sec6}. \textit {We also obtain in Section
\ref{sec6} that if the Livsic characteristic function $W(z)$ of a
bounded dissipative operator with rank one imaginary part satisfies
the condition $W(-z)=W^{-1}(z)$ in some neighborhood of infinity,
then the corresponding Jacobi matrix possesses the property $\RE
b_1=b_2=\cdots =0.$ In Section \ref{sec7} we find the Jacobi matrix
corresponding to the integration operator and establish that the
perturbation of its upper left entry $b_1$ leads to a Jacobi matrix
with a complete system of eigensubspaces.}

We will use the following notation: For a bounded linear operator
$A$ in $\cH$,  $A_R$ and $A_I$ denote its Hermitian components
\[
A_R=\frac{A+A^*}{2}\,,\; A_I=\frac{A-A^*}{2i},\
\]
and $\rho(A)$ and $\sigma(A)$ denote the resolvent set and the
spectrum of $A$.
 By $\dC_+
(resp.,\dC_-)$ we denote the open upper (resp., lower) half-plane.
For a continued fraction
\[
 \frac{c_1}{d_1+\displaystyle\frac{\mathstrut
c_2}{d_2+\displaystyle\frac{\mathstrut c_3}{d_3+\cdots}}}
\]
we use the notation (see \cite{JTh}, \cite{KK}, \cite{Wall})
\[
\frac{c_1}{d_1}\;\raisebox{-3mm}{{\rm
+}}\;\frac{c_2}{d_2}\;\raisebox{-3mm}{{\rm
+}}\;\frac{c_3}{d_3}\;\raisebox{-3mm}{{\rm +}\,\dots}
\]
 {\sl Finally, we would like to mention that the
fundamental research provided by Yuri Berezanskii, Harry Hochstadt,
Fritz Gesztesy and Barry Simon in the theory of self-adjoint Jacobi
matrices as well as by Moshe Livsic in the theory of
non-self-adjoint operators inspired and encouraged us to make a new
step and to develop the direct and inverse spectral analysis of
finite and semi-infinite non-self-adjoint Jacobi and tri-diagonal
matrices with a rank one imaginary part as presented in this paper.}

\section{Bounded prime dissipative operators and their characteristic functions}
\label{sec2}
\subsection{Bounded prime dissipative operators} Let $\cH$ be a complex separable Hilbert space
with the inner product $(\cdot,\cdot)$ and the norm $||\cdot||$. The
operator $A$ is called \textit{dissipative} if $A_I\ge 0$. The
resolvent set $\rho(A)$ of dissipative operator $A$ contains the
open lower half-plane and the estimate
\[
||(A-zI)^{-1}||\le\frac{1}{|\IM z|}, \; \IM z<0
\]
holds. Let
\begin{equation}
\label{ssub}
 \cH_s= \overline{{\rm span}\left\{A^nA_I \cH,
n=0,1,\ldots\right\}}.
\end{equation}
The subspace $\cH_s$ reduces $A$ and $A\uphar(\cH\ominus\cH_s)$ is a
self-adjoint operator. The operator $A\uphar\cH_s$ is called the
\textit{prime part} of $A$.
\begin{definition}
\label{Prime} A bounded linear operator $A$ in a separable Hilbert
space is called prime or completely non-self-adjoint if there is no
reducing invariant subspace on which the operator $A$ is
self-adjoint.
\end{definition}
It is well known \cite{Br} that an operator $A$ is prime in $\cH$ if
and only if
\begin{equation}
\label{prime}   \overline{{\rm span}}\left\{A^nA_I \cH,
n=0,1,\ldots\right\}=\cH.
\end{equation}
\begin{proposition}
\label{SS} Let $A$ be a dissipative operator with a rank one
imaginary part and let $g$ be a vector in $\cH$ such that
\begin{equation}
\label{IMPART} 2A_Ih=(h,g)g,\; h\in\cH.
\end{equation}
Then $A$ is prime if and only if the vector $g$ is cyclic for the
real part $A_R$.
\end{proposition}
\begin{proof}
Suppose that $A$ is prime. Then \eqref{prime} holds.
 Let us prove that $g$ is a cyclic vector for $A_R$. Let
\[
\cH'= \overline{{\rm span}}\left\{(A_R)^n g, n=0,1,\ldots\right\}\ne
\cH.
\]
Then $\cH'$ and $\cH''=\cH\ominus\cH'$ are invariant with respect to
$A_R$. Since $\cH''\subset \Ker A_I$, it follows that
\[
A_R\uphar\cH''=A\uphar\cH''=A^*\uphar\cH'',
\]
 $A^{*n}\cH''\subset\cH$ and $A_IA^{*n}\uphar\cH''=0$ for all
 $n=0,1,\ldots.$ Now from \eqref{prime} we obtain that
 $\cH''=\{0\},$ i.e. $g$ is a cyclic vector for $A_R$.

 Conversely, suppose that the vector $g$ is cyclic for $A_R$, i.e. $\cH'=\cH$. Let
 the subspace $\cH_s$ be defined by \eqref{ssub}. Then
 $$A\uphar(\cH\ominus\cH_s)=A_R\uphar(\cH\ominus\cH_s)$$
  and
 $\cH\ominus\cH_s$ as well as $\cH_s$ reduces $A_R$. Because
 $g\in\cH _s$, we get that $\cH'\subset\cH_s.$ It follows that
 $\cH_s=\cH$, i.e. $A$ is a prime operator.
\end{proof}
A prime dissipative operator has no real eigenvalues and its
non-real eigenvalues belong to the open upper half-plane. It is
known (see \cite{Br}) that the non-real spectrum of an operator with
compact imaginary part consists of eigenvalues of finite algebraic
multiplicities (dimensions of the corresponding root subspaces) and
the limit points of non-real spectrum belong to the spectrum of the
real part of operator. \textit{Note that any eigensubspace of an
operator with a rank one imaginary part which corresponds to a
non-real eigenvalue has a dimension equal to one}. Observe also that
if $\dim\cH=n$ and $A$ is a dissipative operator with a rank one
imaginary part $A_I$, and if \eqref{IMPART} holds, then as it
follows from Definition \ref{prime} and Theorem \ref{Livsic} the
following conditions are equivalent:
\begin{enumerate}
\item $A$ is a prime operator;
\item $A$ has no real eigenvalues;
\item the vectors $g,Ag,\ldots,A^{n-1}g$ are linearly independent.
\end{enumerate}
The following completeness criterion has been established by
M.~S.~Livsic \cite{L}.
\begin{theorem} \cite{L}
\label{Livsic1}. Let $A$ be a bounded prime dissipative operator in
a separable Hilbert space $\cH$ with the finite trace ${\rm Sp}A_I$.
Then the closure of the linear span of all root subspaces coincides
with $\cH$ if and only if
\begin{equation}
\label{trace} \sum\limits_{n}\IM z_n ={\rm Sp} A_I,
\end{equation}
where $\{z_n\}$ is the set of all non-real eigenvalues of $A$
counting multiplicity.
\end{theorem}
\subsection{The Livsic characteristic
function of non-self-adjoint bounded operators} In this subsection
we present basic facts of the Livsic characteristic functions theory
 \cite{Br}, \cite{BrL}, \cite{L}, \cite{LYan}. We restrict
ourselves just to the case of a bounded dissipative operator $A$
with a rank one imaginary part $A_I$.

A vector $g$ satisfying \eqref{IMPART} is called a channel vector
for $A$. The function
\begin{equation}
\label{Char} W(z)=1-i\left((A-zI)^{-1}g,g\right), \; z\in\rho(A),
\end{equation}
is called a \textit{characteristic function of }$A$.  Because an
operator $A$ is bounded and  dissipative the characteristic function
$W(z)$ possesses the following properties:
\begin{enumerate}
\item $\left|W(z)\right|\ge 1$ for $z\in\rho(A)$, $\IM z>0$,
\item $\left|W(z)\right|\le 1$, $\IM z<0$,
\item $\left|W(z)\right|= 1$ for $z\in\rho(A)$, $\IM z=0$,
\item
$ W(z)=1+i \sum_{n=0}^\infty\frac{(A^n g,g)}{z^{n+1}}$ in a
neighborhood of infinity.
\end{enumerate}
In addition the function
\begin{equation}
\label{kernel}
G(z,\xi)=\frac{1-W(z)\ovl{W(\xi)}}{i(z-\bar\xi)}\;\quad\mbox{is a
nonnegative kernel},
\end{equation}
i.e.
\[
 \sum\limits_{i,j=1}^n \lambda_i\overline{\lambda_j}G(z_j,z_i)
 \ge 0, (n=1,\dots)
\]
for any set of points $\{z_i\}_{i=1}^n\subset\rho(A)$ and any vector
$(\lambda_1,\cdots,\lambda_n)\in\dC^n.$
\begin{theorem} \cite{L}.
\label{Livsic2} Let $A_1$ and $A_2$ be bounded dissipative operators
with a rank one imaginary part. If the corresponding characteristic
functions $W_1(z)$ and $W_2(z)$ coincide in some neighborhood of
infinity, then the prime parts of $A_1$ and $A_2$ are unitarily
equivalent.
\end{theorem}
\begin{theorem}\cite{Br}.
\label{Livsic3} Let a function $W(z)$ be holomorphic outside some
closed bounded subset $\cD$ of the closed upper half-plane and
possesses the properties
\begin{enumerate}
\def\labelenumi{\rm (\roman{enumi})}
\item $\left|W(z)\right|\le 1$, $\IM z<0$,
\item $\left|W(z)\right|= 1$ for  $z\in\dC\setminus\cD$, $\IM z=0$,
\item $\lim\limits_{z\to\infty}W(z)=1$.
\end{enumerate}
Then there exists a prime bounded dissipative operator with a rank
one imaginary part whose characteristic function coincides with
$W(z)$.
\end{theorem}

\begin{theorem}
\label{Livsic}\cite{L}. The characteristic function $W(z)$ of a
bounded dissipative operator $A$ with a rank one imaginary part
admits the following multiplicative representation:
\begin{equation}
\label{multrep} W(z)=\prod\limits_{k=1}^N \left(\frac{z-{\bar
z}_k}{z-z_k}\right)^{n_k}\, \cdot\exp{\left(-i\int\limits_{0}^l
\frac{dx}{\alpha(x)-z}\right)},
\end{equation}
where $\{z_k\}_{k=1}^N$ ($N\le\infty$) are the distinct non-real
eigenvalues of $A,$
 $\{n_k\}$ are their algebraic multiplicities,
$\alpha(x)$ is a bounded nondecreasing and continuous from the right
function on the interval $[0,l]$.
\end{theorem}
Observe, that for any operator $A$ and his channel vector $g$
\eqref{IMPART}
\[
\sum\limits_{k=1}^N n_k\IM z_k\le \frac{||g||^2}{2}.
\]
 Moreover, from Theorems \ref{Livsic1} and \ref{Livsic} it follows that
 the closure of the linear span of all root subspaces of operator $A$ coincides with
$\cH$ if and only if
\[
\sum\limits_{k=1}^N n_k\IM z_k=\frac{||g||^2}{2}.
\]
\subsection{Triangular models for bounded dissipative operators with
rank one imaginary parts}\label{TRMOD}
 The next theorem provides a
triangular model for a finite-dimensional prime dissipative operator
with one dimensional imaginary part.
\begin{theorem}\cite{GL}, \cite{L}.
\label{model} Let $\dim \cH=n$ and let $A$  be a dissipative
operator in $\cH$ with a rank one imaginary part $A_I$. Then there
exists an orthonormal basis in $\cH$ in which the matrix of the
operator $A$ takes the form
\[
\dot{A}= \begin{pmatrix} \alpha_1+\frac{i}{2}\beta^2_1 &
i\beta_1\beta_2 & \cdot &
\cdot & \cdot & i\beta_1\beta_n \\
0 & \alpha_2+\frac{i}{2}\beta^2_2 &  \cdot &
\cdot & \cdot & i\beta_2\beta_n \\
\cdot & \cdot & \cdot & \cdot &
\cdot & \cdot  \\
0 & 0 & \cdot & \cdot & \cdot & \alpha_n+\frac{i}{2}\beta^2_n
\end{pmatrix},
\]
where $\{\alpha_k\}_{k=1}^n$ are real numbers and
$\{\beta_k\}_{k=1}^n$ are nonnegative numbers. The operator $A$ is
prime if and only if the numbers $\{\beta_k\}$ are positive. The
characteristic function of the operator $A$ takes the form
\[
W(z)=\prod\limits_{\beta_k\ne
0}\frac{z-\alpha_k+i\beta^2_k/2}{z-\alpha_k-i\beta^2_k/2}.
\]
\end{theorem}
Now consider the case of infinite-dimensional operator.

1) Suppose that a bounded prime dissipative operator $A$ with rank
one imaginary part has a complete system of root subspaces. Then the
characteristic function of $A$ takes the form
\[
\prod\limits_{k=1}^\infty\frac{z-\bar{z}_k}{z-z_k},
\]
where $\{z_k\}$ are the set of not necessarily distinct complex
numbers counting multiplicity and  such that
\begin{enumerate}
\item
$\IM z_k>0$,
\item
 $\sum_{k=1}^\infty\IM z_k<\infty$,
\item
$|z_k|<C.$
\end{enumerate}
From the Livsic approach \cite{L} such an operator $A$ is unitarily
equivalent to the operator $\dot{A}$ in the Hilbert space $l_2(\dN)$
given by the triangular matrix
\[
\dot{A}=
\begin{pmatrix} \alpha_1+\frac{i}{2}\beta^2_1 & i\beta_1\beta_2 &
\cdot &
\cdot & \cdot & i\beta_1\beta_k&\cdot&\cdot&\cdot\\
0 & \alpha_2+\frac{i}{2}\beta^2_2 &  \cdot & \cdot & \cdot &
i\beta_2\beta_k&\cdot&\cdot&\cdot \\ \cdot & \cdot & \cdot & \cdot &
\cdot & \cdot& \cdot & \cdot & \cdot  \\
0 & 0 & \cdot & \cdot & \cdot & \alpha_k+\frac{i}{2}\beta^2_k& \cdot & \cdot & \cdot\\
\cdot&\cdot&\cdot&\cdot& \cdot & \cdot & \cdot& \cdot & \cdot
\end{pmatrix},
\]
where $\{\alpha_k\}_{k=1}^n$ are real numbers and
$\{\beta_k\}_{k=1}^n$ are positive numbers.

2) Suppose that $A$ has only real spectrum. Then the characteristic
function of $A$ takes the form
\[
W(z)=\exp{\left(-i\int\limits_{0}^l \frac{dx}{\alpha(x)-z}\right)}.
\]
Consider the following operator $\dot{B}$ in the Hilbert space
$L_2[0,l]$:
\[
(\dot{B} f)(x)=\alpha(x)f(x)+i\int\limits_{x}^l f(t) dt.
\]
The operator $A$ is unitarily equivalent to the prime part of the
operator $\dot{B}$ \cite{L}.

 3) Let A be an arbitrary bounded prime
dissipative operator with a rank one imaginary part acting on
infinite-dimensional Hilbert space $\cH$. Then its characteristic
function admits  multiplicative representation \eqref{multrep}.
Consider Hilbert space $\sH=\dC^n\oplus L_2[0,l]$ in case when the
non-real spectrum of $A$ is finite and consists of $n$ complex
numbers counting their algebraic multiplicity and the Hilbert space
$\sH=l_2(\dN)\oplus L_2[0,l]$ if the non-real spectrum of $A$ is
infinite. Let $\dot{G} $ be an operator in $\sH$ given by the block
operator matrix
\[
\dot{G}=\begin{pmatrix} \dot{A}&\Gamma\cr 0&\dot{B} \end{pmatrix},
\]
where $\Gamma:L_2[0,l]\to \dC^n $ or $\Gamma:L_2[0,l]\to l_2(\dN)$
is given by the formula
\[
(\Gamma f)(x)=i\,\left(\sqrt{2\,\sum\limits_k\IM
z_k}\,\int\limits_{0}^l f(t)dt\right) \delta_1,
\]
where $\delta_1\in \dC^n$ or $\delta_1\in l_2(\dN)$ is the vector
with the first component equal to 1 and remaining components equal
to zero. Then $A$ is unitarily equivalent to the prime part of
$\dot{G}$ \cite{L}.
\subsection{Linear-fractional transformation of the characteristic
function} Let $T$ be a bounded linear operator acting on some
Hilbert space $\cH$ and let $h$ be a nonzero vector in $\cH$.
Consider the family of rank one perturbations of $T$ given by
\[
T_t=T+t(\cdot,h)\,h,
\]
where $t$ is a complex number. One can easily derive the following
relation
\begin{equation}
\label{RESFOR}
\left((T_t-zI)^{-1}h,h\right)=\frac{\left((T-zI)^{-1}h,h\right)}{1+t\left((T-zI)^{-1}h,h\right)}\,,\;
z\in\rho(T_t)\cap \rho(T).
\end{equation}
Moreover, the number $z_0\in\rho(T)$ is the eigenvalue of $T_t$ if
and only if $z_0$ satisfies the equation
\[
1+t\left((T-z_0I)^{-1}h,h\right)=0.
\]
Let $A$  be a dissipative operator in $\cH$ with a rank one
imaginary part $A_I$ and let $g\in \cH$ such that \eqref{IMPART}
holds.
 Define the
following function
\begin{equation}
\label{V} V(z)=\frac{1}{2}\left((A_R-zI)^{-1}g,g\right),\;
z\in\rho(A_R).
\end{equation}
The function $V(z)$ is holomorphic on the domain $C\setminus [a,b]$,
where $a$ and $b$ are the lower and upper bounds of the spectrum of
self-adjoint operator $A_R$, and posseses the properties
\begin{enumerate}
\item $\displaystyle\frac{V(z)-V^*(z)}{z-\overline{z}}\ge 0,\; z\in\dC_+\cup\dC_-,$
\item $V^*(\overline z)=V(z),$
\item $V(z)=-\frac{1}{2}\sum\limits_{n=0}^\infty\displaystyle\frac{\left((A_R)^n
g,g\right)}{z^{n+1}}$ in a neighborhood of infinity.
\end{enumerate}
The properties (1) and (2) mean that $V(z)$ is a
Herglotz--Nevanlinna function and admits the integral representation
\[
V(z)=\int\limits_{a}^b\frac{d\Omega(t)}{t-z},
\]
where $\Omega(t)=(E(t)g,g)/2\;$ and $E(t)$ is a resolution of
identity for the operator $A_R$. Since
\[
A_R=A-\frac{i}{2}(\cdot,g)\,g,
\]
from \eqref{RESFOR}  and \eqref{Char} it follows that (\cite{Br})
\begin{equation}
\label{FLT} \left\{
\begin{split}
&V(z)=i\,\frac{W(z)-1}{W(z)+1}\\
&W(z)=\frac{1-iV(z)}{1+iV(z)}
\end{split}
\right.,\;
 z\in\rho(A)\cap\rho(A_R).
\end{equation}
The non-real spectrum of $A$ consists of all $z$, $\IM z>0$ which
are solutions of the equation
\[
V(z)=i.
\]
\subsection{The Weyl function of a bounded dissipative operator}
\label{subsec2} In this subsection we define and study the
\textit{Weyl function} of a prime dissipative operator with rank one
imaginary part. In the classical case of the Weyl functions of
self-adjoint differential and difference operators we refer to
\cite{Akh}, \cite{Atk}, \cite{Ber}, \cite{Ber1}, \cite{GS},
\cite{GS1}, \cite{Tes}.
\begin{definition} \label{WFun} Let $A$ be a prime
dissipative operator with a rank one imaginary part. Let $e\in \ran
(A-A^*)$, $||e||=1$. The functions
\begin{equation}
\label{wfuncA}
 m_A(z)=\left((A-zI)^{-1}e,e\right),\; z\in\rho(A)
\end{equation}
and
\begin{equation}
\label{wfuncAr}
 m_{A_R}(z)=\left((A_R-zI)^{-1}e,e\right),\; z\in\rho(A_R)
\end{equation}
are said to be the Weyl functions of the operator $A$ and $A_R$,
correspondingly.
\end{definition}
 Since $e\in\ran(A-A^*)$ and $||e||=1$, the operator $A_I$
takes the form
\[
A_If=l(f,e)e,
\]
where $l>0$. Observe that by Proposition \ref{SS} the vector $e$ is
 cyclic for $A$ and $A_R$. The function $m_A(z)$ is holomorphic on
the resolvent set of $A$ which includes the open lower half-plane
and the neighborhood of infinity of the form $|z|>||A||$ and in this
neighborhood
\[
m_A(z)=-\sum\limits_{k=0}^\infty\frac{(A^ke,e)}{z^{k+1}}.
\]
The function $m_{A_R}(z)$ is a Herglotz-Nevanlinna function
holomorphic on the resolvent set of $A_R$ and for $|z|>||A_R||$
\[
m_{A_R}(z)=-\sum\limits_{k=0}^\infty\frac{((A_R)^ke,e)}{z^{k+1}}.
\]
 Relations \eqref{IMPART} and \eqref{V}
yield
\[
m_{A_R}(z)=\frac{1}{l}\, V(z), \; z\in\rho(A_R).
\]
 From \eqref{IMPART} and \eqref{Char} it follows that
the characteristic function $W(z)$ of $A$ and the Weyl function
$m_A(z)$ are connected via relations
\begin{equation}
\label{CONN1}
 W(z)=1-2il\,m_A(z),\;m_A(z)=\frac{i}{2l}(W(z)-1),\; z\in\rho(A).
\end{equation}
The relation \eqref{RESFOR} gives the following connections between
the Weyl functions of $A$ and $A_R$
\begin{equation}
\label{CONNECT}
\begin{split}
&m_A(z)=\frac{i\,m_{A_R}(z)}{i-l\, m_{A_R}(z)},\\
&m_{A_R}(z)=\frac{i\,m_{A}(z)}{i+l\, m_{A}(z)},\;
z\in\rho(A)\cap\rho(A_R).
\end{split}
\end{equation}
From \eqref{CONNECT} it follows that the non-real eigenvalues of $A$
coincide with poles of the function $m_A(z)$ and zeroes of the
function $i-l\,m_{A_R}(z)$. The following theorem directly follows
from Theorems \ref{Livsic2} and relations \eqref{CONN1}.
\begin{theorem}
\label{UNEQ} Let $A_1$ and $A_2$ be two prime dissipative operators
with rank one imaginary parts. If the Weyl functions of $A_1$ and
$A_2$ coincide in some neighborhood of infinity, then the operators
$A_1$ and $A_2$ are unitarily equivalent.
\end{theorem}
 Suppose that a Hilbert space $\cH$ is decomposed
as $\cH=\cH_1\oplus\cH_2$. Then every bounded operator $A$ in $\cH$
has the block-operator matrix representation with respect to this
decomposition
\[
A=\begin{pmatrix}A_{11}&A_{12}\cr A_{21}&A_{22}\end{pmatrix}.
\]
According to the Schur-Frobenius formula the relation
\begin{equation}
\label{ShF} P_{\cH_{1}}(A-zI)^{-1}\uphar\cH_1=\left(-zI+
A_{11}-A_{12}(A_{22}-zI)^{-1}A_{21}\right)^{-1}
\end{equation}
holds for $z\in\rho(A)\cap\rho(A_{22}).$ Here $P_{\cH_{1}}$ is the
orthogonal projection in $\cH$ onto $\cH_1$. Let $A$ be a
dissipative operator with a rank one imaginary part and let
$\cH_1=\ran (A-A^*)$. Then $\cH_2=\Ker(A-A^*)$, the operator
$A_{22}$ is self-adjoint in $\cH_2$, and $A_{12}=A^*_{21}.$ It
follows that
\[
A_{11}e=b e,\; A_{12}f=a(f,h)e,\; f\in \cH_2,\; A_{21}e=\bar a h,
\]
where $e\in \ran(A-A^*)$ and $h\in\Ker(A-A^*)$ are  the unit vectors
and $\IM b>0.$  Hence \eqref{ShF} takes the form for
$z\in\rho(A)\cap\rho(A_{22})$
\begin{equation}
\label{SSSS}
\left((A-zI)^{-1}e,e\right)=\frac{1}{-z+b-|a|^2\left((A_{22}-zI)^{-1}h,h\right)}.
\end{equation}
 If $A$ is a prime operator, then $a\ne 0$ and in this case for $z\in\rho(A)\cap\rho(A_{22})$
\begin{equation}
\label{1000}
\left((A_{22}-zI)^{-1}h,h\right)=\frac{1}{|a|^2}\left(-z+b-\frac{1}{\left((A-zI)^{-1}e,e\right)}\right).
\end{equation}
It follows that
\[
\begin{split}
&\left((A_{22}-zI)^{-1}h,h\right)=\frac{1}{|a|^2}
\left(b-z-\frac{1}{m_A(z)}\right),\; z\in
\rho(A)\cap\rho(A_{22}),\\
&\left((A_{22}-zI)^{-1}h,h\right)=\frac{1}{|a|^2}\left(\RE b -
z-\frac{1}{m_{A_R}(z)}\right),\;z\in\rho(A_R)\cap\rho(A_{22}).
\end{split}
\]
In addition, the number $z_0$, $\IM z_0>0$ is the eigenvalue of $A$
if and only if $z_0$ satisfies the equation
\[
b-z-|a|^2\left((A_{22}-zI)^{-1}h,h\right)=0.
\]
The corresponding eigenspace is
\[
\span\left\{e-\overline a (A_{22}-z_0I)^{-1}h\right\}.
\]
 The next theorem gives the analytical characterization
of the Weyl functions of prime dissipative operators with a rank one
imaginary part.
\begin{theorem}
\label{WeylFunction} Let $M(z)$ be a function holomorphic outside
some closed bounded subset $\cD$ of the closed upper half-plane. The
following statements are equivalent:
\begin{enumerate}
\def\labelenumi{\rm (\roman{enumi})}
\item The function $M(z)$ is the Weyl function of some bounded prime dissipative operator
with a rank one imaginary part.
\item
\begin{enumerate}
\item The function $M(z)$ has the asymptotic expansion at infinity
\begin{equation}
\label{ASYMPTOT}
M(z)=-\frac{1}{z}-\frac{b}{z^2}+O\left(\frac{1}{z^3}\right),
\end{equation}
where $\IM b>0$.
\item The function
\[
W(z)=1-2i\,\IM b\, M(z)
\]
is the characteristic function of some dissipative bounded operator
with a rank one imaginary part.
\end{enumerate}
\item
\begin{enumerate}
\item The function $M(z)$ has the expansion \eqref{ASYMPTOT} at
infinity.
\item The function
\begin{equation}
\label{KERN}
 \cK(z,\xi):=\frac{M(z)-\ovl{M(\xi)}+2i\,\IM b\,
M(z)\ovl{M(\xi)}}{z-\ovl\xi}
\end{equation}
is a nonnegative kernel.
\end{enumerate}
\item
\begin{enumerate}
\item The function $M(z)$ has the expansion \eqref{ASYMPTOT} at
infinity.
\item The function
\[
Q(z):=\frac{i\,M(z)}{i+\IM b\, M(z)}
\]
has analytic continuation onto the outside of some bounded interval
of the real axis as the Herglotz--Nevanlinna function .
\end{enumerate}
\item
\begin{enumerate}
\item The function $M(z)$ has the expansion \eqref{ASYMPTOT} at
infinity.
\item The function
\[
\sM(z):=-\frac{1}{M(z)}+b
\]
has analytic continuation onto the outside of some bounded interval
of the real axis as the Herglotz--Nevanlinna function.
\end{enumerate}
\end{enumerate}
\end{theorem}
\begin{proof}
Let $M(z)$ be the Weyl function of some bounded prime dissipative
operator $A$ with a rank one imaginary part. From \eqref{wfuncA} we
get
\[
\begin{split}
&M(z)=\left((A-zI)^{-1}e,e)\right)=-\sum\limits_{k=0}^\infty\frac{(A^ke,e)}{z^{k+1}}\\
&=\frac{1}{z}-\frac{b}{z^2}+O\left(\frac{1}{z^3}\right),\;
|z|>||A||,
\end{split}
\]
where $e\in\ran(A-A^*),\; ||e||=1,$ $b=(Ae,e)$, and $\IM b>0$. So,
\eqref{ASYMPTOT} holds. Since $A_Ie=\IM b \,e$, from \eqref{Char} it
follows that the function $W(z)=1-2i\,\IM b\, M(z)$ coincides with
the characteristic function of $A$. Hence
\[
 \cK(z,\xi)=\frac{M(z)-\ovl{M(\xi)}+\IM b\, M(z)\ovl{M(\xi)}}{z-\ovl\xi}=
\frac{1-W(z)\overline{W(\xi)}}{i(z-\bar\xi)}
\]
is a nonnegative kernel. By \eqref{CONNECT} the function
\[
Q(z)=\frac{i\,M(z)}{i+\IM b\, M(z)}, \; z\in\rho(A)\cap\rho(A_R)
\]
coincides with the Weyl function of $A_R$ and therefore, has
analytic continuation onto the outside of some bounded interval of
the real axis as the Herglotz--Nevanlinna function. By \eqref{1000},
 the function
 \[
\sM(z)=-\frac{1}{M(z)}+b
\]
also has analytic continuation onto the outside of some bounded
interval of the real axis as the Herglotz--Nevanlinna function.
Thus, the statement (i) implies the statements (ii), (iii), (iv),
and (v).

Now suppose that the statement (ii) holds true, i.e. there exist a
Hilbert space $\cH$ and a bounded prime dissipative operator $A$ in
$\cH$ with a rank one imaginary part such that
\[
1-i ((A-zI)^{-1}g,g)=1-2i\IM bM(z),
\]
where $g\in\ran A_I$ and $A_I=(\cdot,g)g$. It follows that
$||g||^2=2\IM b$ and
\[
M(z)=((A-zI)^{-1}e,e),\quad\mbox{where}\quad e=\frac{g}{\sqrt{2\IM
b}}\,.
\]
By definition \ref{WFun} the function $M(z)$ is the Weyl function of
$A$. So, (ii)$\Rightarrow$ (i). If (iii) holds, then for
$W(z)=1-2i\IM b M(z)$ from \eqref{KERN} we obtain that
\[
\cK(z,\xi)=\frac{1-W(z)\overline{W(\xi)}}{i(z-\bar\xi)}.
\]
Since $\cK(z,\xi)$ is a nonnegative kernel, we get that
$\left|W(z)\right|\le 1$, $\IM z<0$, and $\left|W(z)\right|= 1$ for
$z\in\dC\setminus\cD$, $\IM z=0$. In addition, from \eqref{ASYMPTOT}
we get that $\lim\limits_{z\to\infty}W(z)=1$. By Theorem
\ref{Livsic3} the function $W(z)$ is the characteristic function of
some bounded prime dissipative operator with rank one imaginary
part. Thus (iii)$\Rightarrow$(ii).

Let us show that (iv)$\Rightarrow$(ii). Define the function
\[
W(z)=\frac{1-i\IM b\, Q(z)}{1+i\IM b\,Q(z)}.
\]
Then $W(z)=1-2i\IM bM(z)$. Since $Q(z)$ is a Herglotz--Nevanlinna
function, we have $\left|W(z)\right|\le 1$, $\IM z<0$,
$\left|W(z)\right|= 1$ for $z\in\dC\setminus\cD$, $\IM z=0$, and
$\lim\limits_{z\to\infty}W(z)=1$. By Theorem \ref{Livsic3} the
function $W(z)$ is the characteristic function of some bounded prime
dissipative operator with a rank one imaginary part.

Suppose that (v) holds true. Because
\[
\sM(z)=-\frac{1}{M(z)}+b
\]
is a Herglotz--Nevanlinna function, the function
\[
\frac{\sM(z)-\overline{\sM(\xi)}}{z-\bar\xi}
\]
is a nonnegative kernel.  From \eqref{KERN} it follows that
\[
\cK(z,\xi)=M(z)\overline{M(\xi)}\,\frac{\sM(z)-\overline{\sM(\xi)}}{z-\bar\xi}.
\]
Therefore, $\cK(z,\xi)$ is a nonnegative kernel and (v)$\Rightarrow$
(ii).
\end{proof}
\section{Jacobi matrices and their Weyl
functions} \label{sec3}
\subsection{Complex Jacobi matrices and corresponding Jacobi operators}
Let $\{b_k\}$, $\{a_k\}$, $a_k\ne 0$, $k\ge 1$ be complex numbers. A
tri-diagonal matrix of the form
\begin{equation}
\label{22} J=\begin{pmatrix} b_1 & a_1 & 0  &0  & 0 & \cdot &
\cdot & \cdot \\
a_1 & b_2 & a_2 & 0 &0& \cdot &
\cdot & \cdot \\
0    & a_2 & b_3 & a_3 &0& \cdot &
\cdot & \cdot \\
\cdot & \cdot & \cdot & \cdot &
\cdot & \cdot & \cdot &\cdot \\
\cdot & \cdot & \cdot & \cdot &
\cdot & \cdot & \cdot &a_{n-1} \\
\cdot & \cdot & \cdot & \cdot &\cdot & 0 & a_{n-1} & b_n
\end{pmatrix}
\end{equation}
is called finite complex Jacobi matrix and a matrix
\begin{equation}
\label{222} J=\begin{pmatrix} b_1 & a_1 & 0 &0   & 0 & \cdot &
\cdot & \cdot \\
a_1 & b_2 & a_2 & 0 &0& \cdot &
\cdot & \cdot \\
0    & a_2 & b_3 & a_3 &0& \cdot &
\cdot & \cdot \\
\cdot & \cdot & \cdot & \cdot & \cdot & \cdot & \cdot &\cdot
\end{pmatrix}
\end{equation}
is called semi-infinite complex Jacobi matrix. The case of real
entries $\{b_k\}$ and positive entries $\{a_k\}$ corresponds to the
classical symmetric Jacobi matrix \cite{Akh}, \cite{S}, \cite{Tes}.
We will call such  matrix a self-adjoint Jacobi matrix. A finite
$n\times n$ Jacobi matrix determines a linear operator (the Jacobi
operator) in the Hilbert space $\dC^n$.  Let $\cC_0$ be the linear
manifold of vectors in $l_2(\dN)$ with finite support. A
semi-infinite Jacobi matrix determines two linear operators in
$l_2(\dN)$ given by the formal matrix product $J\cdot x$. The first
operator is defined in $\cC_0$. This operator is densely defined and
is closable. Let $[J]_{min}$ be its closure. The second operator
$[J]_{max}$ has the domain
\[
\dom([J]_{max})=\{x\in l_2(\dN): J\cdot x\in l_2(\dN)\}.
\]
A semi-infinite Jacobi matrix is called proper if
$[J]_{min}=[J]_{max}$ (cf. \cite{Bec1}).
 It is well known (cf.
\cite{AkhGl}, \cite{Ber}) that $\dom[J]_{min}=l_2(\dN)$ and
$[J]_{min}$ is bounded in $l_2(\dN)$ if and only if the entries
$\{a_k\}$ and $\{b_k\}$ are uniformly bounded, i.e. condition
\eqref{Bound} is fulfilled. The tri-diagonal matrix is compact if
and only if
\[
\lim\limits_{k\to\infty}b_k=\lim\limits_{k\to\infty}a_k=0.
\]
A Jacobi matrix is called bounded if \eqref{Bound} holds.

Because \begin{equation} \label{coord}
\begin{split}
&\left(J^k\delta_1\right)_{k+1}=a_k a_{k-1}\ldots a_1,\\
&\left(J^k\delta_1\right)_m=0,\; m\ge k+2,
\end{split}
\end{equation}
and $a_k\ne 0$, the vectors $\delta_1, J\delta_1,\ldots,
J^k\delta_1,\ldots$ are linearly independent and moreover the vector
$\delta_1$ is cyclic for the operator $[J]_{min}$ in $l_2(\dN)$.

The system of second order difference equations
\begin{equation}
\label{SYSTEM}
\begin{split}
& a_k P_{k+1}(z)+b_kP_k(z)+a_{k-1}P_{k-1}(z)=zP_k(z),\; k\ge 1,\;
a_{0}:=1,\\
&  \;P_0(z)=0, \;P_1(z)=1
\end{split}
\end{equation}
determines polynomials $P_k(z)$ of degree $k-1$ (in the case of
$n\times n$ Jacobi matrix we define $a_n:=1$). Note that
\begin{equation}
\label{POLYNOM} P_{k+1}(z)=\frac{1}{a_1\cdots
a_{k}}z^{k}+\quad\mbox{lower degree in}\: z,\; k=1,\ldots, n.
\end{equation}
\subsection{Finite complex Jacobi matrices}
 The next proposition probably is well-known. We give a proof for
 completeness.
\begin{proposition}
\label{eigcores} Let $J$ be an $n\times n$ complex Jacobi matrix
with $a_k\ne 0$ for all $k$ and let polynomial
$\{P_k(z)\}_{0}^{n+1}$ be defined by \eqref{SYSTEM}. Then $z_0$ is
an eigenvalue of $J$ if and only if $P_{n+1}(z_0)=0$. Moreover, if
$z_0$ is a root of the polynomial $P_{n+1}(z)$ of the multiplicity
$l$, then vectors in $\dC^n$ defined as follows
\begin{equation}
\label{corresp}
\begin{split}
&e_0=\begin{pmatrix} 1\cr P_2(z_0)\cr\vdots\cr\vdots\cr P_{n}(z_0)
\end{pmatrix},\; e_1=\begin{pmatrix} 0\cr P'_2(z_0)\cr\vdots\cr\vdots\cr P'_{n}(z_0)
\end{pmatrix},\;e_2=\frac{1}{2!}\begin{pmatrix} 0\cr 0\cr P''_{3}(z_0)\cr\vdots\cr P''_{n}(z_0)
\end{pmatrix},\;\ldots,\\
& e_{l-1}=\frac{1}{(l-1)!}\begin{pmatrix} 0\cr 0\cr\vdots\cr 0\cr
P^{(l-1)}_{l}(z_0)\cr\vdots \cr P^{(l-1)}_{n}(z_0) \end{pmatrix}
\end{split}
\end{equation}
satisfy the relations
\[
(J-z_0I)e_0=0,\; (J-z_0I)e_{k+1}=e_k,\; k=0,1,\ldots, l-2.
\]
In addition, the relation
\begin{equation}
\label{CHARPOL} P_{k+1}(z)=a_1\, a_2\,\cdots a_k\, {\rm
det}\left(zI-J_{[1,k]}\right), \; k=1,\ldots, n
\end{equation}
holds, where $J_{[1,k]}$ is the $k\times k$ upper left corner of
$J$.
\end{proposition}
\begin{proof} Let $z_0$ be an eigenvalue of $J$ and let
$\begin{pmatrix}y_1&y_2&\ldots&y_n\end{pmatrix}^T$ be the
corresponding eigenvector. Then we have the system of linear
equations
\[
\left\{
\begin{split}
&b_1y_1+a_1 y_2=z_0y_1\\
&a_1y_1+b_2 y_2+a_2y_3=z_0y_2\\
&\ldots\dots\ldots\ldots\ldots\ldots\ldots\ldots\\
&a_{n-2}y_{n-2}+b_{n-1}y_{n-1}+a_{n-1}y_n=z_0y_{n-1}\\
&a_{n-1}y_{n-1}+b_ny_n=z_0y_n .
\end{split}
\right.
\]
Since $a_k\ne 0$ for all $k=1,2\ldots, n-1$, we can express $y_2$,
$y_3$, $\ldots$, $y_n$ linearly through $y_1$. It follows that
$y_1\ne 0$ and we may put $y_1=1$. Comparing with \eqref{SYSTEM} we
obtain that $y_2=P_2(z_0)$, $\ldots,$ $y_n=P_{n}(z_0)$, and
$P_{n+1}(z_0)=0.$

Conversely, if $P_{n+1}(z_0)=0$, then from \eqref{SYSTEM} it follows
that $z_0$ is an eigenvalue of $J$ and the vector
$e_0=(1,P_2(z_0),\ldots, P_n(z_0))^T$ is corresponding eigenvector.

Let $z_0$ be a root of $P_{n+1}(z)$ of the multiplicity $l$. Then
$$P_{n+1}(z_0)=P'_{n+1}(z_0)=\ldots=P^{(l-1)}_{n+1}(z_0)=0,\;
P^{(l)}_{n+1}(z_0)\ne 0.$$ Differentiating \eqref{SYSTEM} we get
\begin{equation}
\label{DER}
 \left\{
\begin{split}
& b_1P'_1(z)+a_1 P'_2(z)=zP'_1(z)+P_1(z)\\
&a_1P'_1(z)+b_2P'_2(z)+a_2P'_3(z)=zP'_2(z)+P_2(z)\\
&\ldots\dots\ldots\ldots\ldots\ldots\ldots\ldots\ldots\dots\ldots\ldots\ldots\ldots\ldots\ldots\\
&a_{n-2}P'_{n-2}(z)+b_{n-1}P'_{n-1}(z)+a_{n-1}P'_{n}(z)=zP'_{n-1}(z)+P_{n-1}(z)\\
&a_{n-1}P'_{n-1}(z)+b_nP'_n(z)+a_{n}P'_{n+1}(z)=zP'_{n}(z)+P_n(z)
\end{split}
\right..
\end{equation}
Recall that $P_1(z)=1$ and $a_n=1$. Substituting $z=z_0$ we obtain
that the vector $e_1=\begin{pmatrix} 0& P'_2(z_0)&\ldots&
P'_{n}(z_0)\end{pmatrix}^T$ satisfies $(J-z_0I)e_1=e_0$. Continuing
the differentiation of \eqref{DER}, we prove that the vectors given
by \eqref{corresp} satisfy the equalities $(J-z_0 I)e_{k+1}=e_k$,
 for $k=1,\ldots, l-2$. This means that the number $z_0$ is the root
of the characteristic polynomial ${\rm det}(J-zI)$ and the
multiplicity of this root is greater or equal $l$.

Suppose now that the number $z_0$ is a root of the characteristic
polynomial ${\rm det}(J-zI)$ of the multiplicity $l>1$. Then as was
proved above the number $z_0$ is the simple eigenvalue of $J$, $z_0$
is a root of the polynomial $P_{n+1}(z)$ and $e_0$ is the
corresponding eigenvector. Suppose that $x=(x_1,x_2,\ldots, x_n)^T$
satisfies the equation $(J-z_0I)x=e_0$. Since the vectors $x+\lambda
e_0$ satisfy the same equation for arbitrary $\lambda$, we may
assume that $x_1=0$. Then we obtain the linear system
\[
 \left\{
\begin{split}
& a_1 x_2=P_1(z_0)\\
&b_2x_2+a_2x_3=zx_2+P_2(z_0)\\
&\ldots\dots\ldots\ldots\ldots\ldots\ldots\ldots\ldots\dots\ldots\ldots\ldots\ldots\ldots\ldots\\
&a_{n-2}x_{n-2}+b_{n-1}x_{n-1}+a_{n-1}x_{n}=zx_{n-1}+P_{n-1}(z_0)\\
&a_{n-1}x_{n-1}+b_nx_n+a_{n}x_{n+1}=zx_{n}+P_n(z_0)
\end{split}
\right..
\]
Comparing the last linear system with \eqref{DER} when $z=z_0$ and
taking into account that $P'_1(z)=0$ for all $z$, we get that
$x_2=P'_2(z_0),\ldots, x_n=P'_n(z_0)$ and $P'_{n+1}(z_0)=0.$ So,
$x=e_1$. Following the same way as above, we get that
$P''_{n+1}(z_0)=\cdots =P^{(l-1)}_{n+1}(z_0)=0,$ i.e. the
multiplicity of $z_0$ as a root of the polynomial $P_{n+1}(z)$ is
greater or equal $l$.

Thus, we proved that the roots and their multiplicities of the
polynomial $P_{n+1}(z)$ coincide with the roots and their
multiplicities of the characteristic polynomial ${\rm det}(J-zI)$.

The same is true for ${\rm det}(J_{[1,k]}-zI)$ and for the
polynomial $P_{k+1}(z)$, $k=0,1,\dots, n$. Taking into account
\eqref{POLYNOM}, we get that \eqref{CHARPOL} holds.
 \end{proof}
\begin{proposition}
\label{COMMON} Let $J$ be an $n\times n$ complex Jacobi matrix
$a_k\ne 0$ for all $k$. Then the matrices $J$ and $J_{[1,n-1]}$ have
no common eigenvalues.
\end{proposition}
\begin{proof}
Suppose that $z_0$ is an eigenvalue of the matrices $J$ and
$J_{[1,n-1]}.$ Then by Proposition \ref{eigcores} we have
$P_{n}(z_0)=P_{n+1}(z_0)=0.$ Since
$$a_{n-1}P_{n-1}(z_0)+b_n
P_n(z_0)+a_{n+1}P_{n+1}(z_0)=z_0 P_n(z_0),$$
 we get $P_{n-1}(z_0)=0$
and from \eqref{SYSTEM} it follows that $P_{n-2}(z_0)=\cdots=
P_{1}(z_0)=0.$ But $P_1(z)=1$ for all $z$. Contradiction. Thus, the
matrices $J$ and $J_{[1,n-1]}$ have no common eigenvalues.
\end{proof}
\subsection{Semi-infinite Jacobi matrices}
In the following we will consider bounded Jacobi semi-infinite
matrices of the form \eqref{222}. For the corresponding Jacobi operator
\begin{equation}
\label{EIGENV}
\begin{array}{l}
\mbox{the number}\;z\in\dC\;\mbox{is an eigenvalue of}\;
 J\iff\\
\qquad\qquad\qquad\sum\limits_{k=1}^\infty\left|P_k(z)\right|^2<\infty
\end{array}
\end{equation}
and the corresponding eigenvectors are $\lambda\left(P_1(z),
P_2(z),\ldots\right)^T$, $\lambda\in\dC$. It was established in
\cite{Bec} that
\[
\rho(J)=\left\{ z\in\dC:\sup\limits_{n\ge
1}\frac{\sum_{k=1}^n\left|P_k(z)\right|^2}
{|a_n|^2\left[\left|P_n(z)\right|^2+\left|P_{n+1}(z)\right|^2\right]}<\infty\right\}.
\]
The function
\begin{equation}
\label{SpFunc} m_J(z)=\left((J-zI)^{-1}\delta_1,\delta_1\right)
\end{equation}
is called the \textit{Weyl function} of $J$.  Note that the
definition of the Weyl function in the form \eqref{SpFunc} for a
non-self-adjoint Jacobi matrix is given in \cite{Bec1}. The Weyl
function has the following Taylor expansion at infinity
\[
m_J(z)=-\sum\limits_{l=0}^\infty\frac{\left(J^l\delta_1,\delta_1\right)}{z^{l+1}}
\]
and the continued fraction ($J$-fractions) expansion \cite{JTh},
\cite{Wall}
\begin{equation}
\label{CF} m_J(z)= \frac{-1}{z-b_1}\;\raisebox{-3mm}{{\rm
+}}\;\frac{-a^2_1}{z-b_2}\;\raisebox{-3mm}{{\rm
+}}\;\frac{-a^2_2}{z-b_3} \;\raisebox{-3mm}{{\rm +}\,\ldots}\;
\raisebox{-3mm}{{\rm +}}\;\frac{-a^2_{n-1}}
{z-b_n}\;\raisebox{-3mm}{{\rm +}\,\dots}.
 \end{equation}
It follows that
\[
m_J(z)\sim-\frac{1}{z}-\frac{b_1}{z^2}-\frac{b^2_1+a^2_1}{z^3}+O\left(\frac{1}{z^4}\right),\;
z\to\infty.
\]
\subsection{ The $m$-function
approach in inverse spectral problems for self-adjoint Jacobi
matrices} \label{subsec3} It is well-known  \cite{H1} that a finite
real $n\times n$ Jacobi matrix $J$ has $n$ distinct real
eigenvalues. A bounded real  semi-infinite Jacobi matrix determines
a bounded self-adjoint Jacobi operator in $l_2(\dN)$ with the simple
spectrum. If $E(t)$ is the orthogonal resolution of identity for the
operator $J$, then $d\rho(t)=d(E(t)\delta_1,\delta_1)$ is a
probability measure supported at $n$ points in the case of an
$n\times n$ Jacobi matrix and on a finite interval of the real axis
in the semi-infinite case. The probability measure $d\rho$ is called
the \textit{spectral measure} of $J$. The corresponding Weyl
function
\[
m_J(z)=\int\frac{d\rho(t)}{t-z}.
\]
belongs to Herglotz-Nevanlinna class.

The inverse spectral problems for finite and semi-infinite
self-adjoint Jacobi matrices were studied in \cite{GS}, \cite{H1},
\cite{H2}, \cite{H3}, \cite{HL}, \cite{Stone}. The following theorem
is established by M.~Stone \cite{Stone} for possibly unbounded
self-adjoint operator.
\begin{theorem} \cite{Akh}, \cite{Stone}.
\label{Stone} Let $A$ be a self-adjoint operator with simple
spectrum in a separable Hilbert space $\cH$. Then there exists an
orthonormal basis in $\cH$ in which the matrix of $A$ is a Jacobi
matrix with the conditions
\begin{equation}
\label{33} b_k \quad\mbox{are real numbers},\; a_k>0 \quad\mbox{for
all} \quad k
\end{equation}
\end{theorem}
 The construction of the corresponding self-adjoint Jacobi
 matrix can be provided by orthogonal normalization of $1, x,
x^2,\ldots$ with respect to the measure $d\rho$ or by means of
$m$-functions (see \cite{GS}).  We briefly describe the later
approach for finite matrix and will keep notations of \cite{GS}.
Denote by $J_{[k,n]}$, $ k=2,\ldots, n$ the Jacobi matrix obtained
from $J$ by deleting $k-1$ top rows and $k-1$ left columns of $J$,
$J_{[1,n]}=J$.
 Let $\delta_k$ be $k\times 1$ column $\delta_k=\begin{pmatrix}0&
\ldots& 0& 1& 0&\ldots&0\end{pmatrix}^T$. Define the functions
\begin{equation} \label{mk}
m_+(z,k-1)=\left((J_{[k,n]}-zI)^{-1}\delta_{k},\delta_{k}\right),\;
k=1,\ldots, n-1.
\end{equation}
Thus $m_+(z,0)=m_J(z)$. As in \cite{GS} we will use the notation
$m_+(z)$ for the function $m_J(z)$.

 Since the functions $m_+(z,k-1)$
has the expansion in a neighborhood of infinity
\[
m_+(z,k-1)=-\sum\limits_{l=0}^\infty\frac{\left(J^l_{[k,n]}\delta_k,\delta_k\right)}{z^{l+1}},
\]
from \eqref{22} it follows that
\begin{equation}
\label{expansion}
\begin{split}
& m_+(z,k-1)\sim-\frac{1}{z}-\frac{b_k}{z^2}-\frac{b^2_k+a^2_k}{z^3}+O\left(\frac{1}{z^4}\right),\;k=1,\ldots, n-1,\\
&m_+(z,n-1)\sim
-\frac{1}{z}-\frac{b_n}{z^2}-\frac{b^2_n}{z^3}+O\left(\frac{1}{z^4}\right),
\end{split}
\end{equation}
and \eqref{SSSS} yields the following relations
\begin{equation}
\label{connect} a^2_k m_+(z,k)+\frac{1}{m_+(z, k-1)}=b_k-z,\;
k=1,\ldots, n.
\end{equation}
\begin{theorem}
 \label{recov} Every $n$-point probability measure is
 the spectral measure of a unique $n\times n$ Jacobi matrix.
\end{theorem}
\begin{proof} We will follow the approach considered in \cite{GS}.
Let $d\rho$ be a probability measure supported at $n$ real points
$t_1, \ldots, t_n$. Define the Herglotz-Nevanlinna function
\[
m(z)=\int
\frac{d\rho(t)}{t-z}=\sum\limits_{k=1}^n\frac{\alpha_k}{t_k-z},
\]
where $\alpha_k>0$ and $\sum\limits_{k=1}^n\alpha_k=1$. Then in the neighborhood of
infinity one has
\[
m(z)\sim-\frac{1}{z}-\frac{b_1}{z^2}-\frac{b^2_1+a^2_1}{z^3}+O\left(\frac{1}{z^4}\right),
\]
where
\[
b_1=\sum\limits_{k=1}^n \alpha_k
t_k,\;b^2_1+a^2_1=\sum\limits_{k=1}^n \alpha_k t^2_k\,.
\]
One can find also $b_1$ and $a_1>0$ by
\[
\begin{split}
&b_1=-\lim\limits_{z\to\infty}z^2\left(m(z)+\frac{1}{z}\right),\\
&a^2_1=-b^2_1-\lim\limits_{z\to\infty}z^3\left(m(z)+\frac{1}{z}+\frac{b^2_1}{z^2}\right).
\end{split}
\]
Define the Herglotz-Nevanlinna function $m_+(z,1)$ by
\[
m_+(z,1)=\frac{1}{a^2_1}\left(b_1-z-\frac{1}{m(z)}\right).
\]
This function has an expansion
\[
m_+(z,1)\sim-\frac{1}{z}-\frac{b_2}{z^2}-\frac{b^2_2+a^2_2}{z^3}+O\left(\frac{1}{z^4}\right)
\]
in the neighborhood of infinity. Find $b_2$ and $a_2>0$ by means of limits at infinity
and define $m_+(z,2)$ similarly. Continuing this process up to
\[
m_+(z,n-1)=\frac{1}{a^2_{n-1}}\left(b_{n-1}-z-\frac{1}{m_+(z,
n-2)}\right)=\frac{1}{b_n-z}
\]
and finding $b_n$ by
\[
b_n=-\lim\limits_{z\to\infty}z^2\left(m_+(z,n-1)+\frac{1}{z}\right),
\]
we obtain the real numbers $\{b_k\}_{k=1}^n$ and the positive
numbers $\{a_k\}_{k=1}^{n-1}$. Let us construct a self-adjoint Jacobi
matrix $J$ of the form \eqref{22}. Then the $m$-function
\[
m_+(z,0)=\left((J-zI)^{-1}\delta_{1},\delta_{1}\right)
\]
coincides with $m(z)$.
\end{proof}
\begin{remark}
\label{CFr} The following  $J$-fraction expansion of the function
$m_+(z)$
\begin{equation}
\label{CF0} m_+(z)= \frac{-1}{z-b_1}\;\raisebox{-3mm}{{\rm
+}}\;\frac{-a^2_1}{z-b_2}\;\raisebox{-3mm}{{\rm
+}}\;\frac{-a^2_2}{z-b_3} \;\raisebox{-3mm}{{\rm +}\,\ldots}\;
\raisebox{-3mm}{{\rm +}}\;\frac{-a^2_{n-1}} {z-b_n}
\end{equation}
holds.
\end{remark}

 For semi-infinite Jacobi matrices with conditions
\eqref{33} and \eqref{Bound} the approach in \cite{GS} is  based on
the representation of $m(z)$ as a continued $J$-fraction of the form
\eqref{CF} and the following result:
\begin{theorem} \cite{GS}.
\label{probmeas} Suppose that $m(z)=\int (t-z)^{-1} d\rho(t)$, where
$d\rho$ is a probability measure on $[-C,C]$ whose support contains
more than one point. Let
\[
b_1=\int d\rho(t), \; a_1^2=\int t^2 d\rho(t)-b^2_1,
\]
and let
\[
m_1(z)=\frac{1}{a^2_1}\left(b_1-z-\frac{1}{m(z)}\right).
\]
Then $m_1(z)=\int (t-z)^{-1}d\rho_1(t)$, where $d\rho_1$ is a
probability measure also supported on $[-C,C]$.
\end{theorem}
We give one more application of the $m$-function approach which we
will use later.
\begin{proposition}
\label{odd} Let $m(z)$ be a Herglotz-Nevanlinna function which is
odd in some neighborhood of infinity and $\lim_{z\to
\infty}zm(z)=-1$. Then the diagonal entries of the corresponding
self-adjoint Jacobi matrix are equal zero.
\end{proposition}
\begin{proof}  Since $m(-z)=-m(z)$, the function $m(z)$ has the following Taylor's
expansion at infinity
\begin{equation}
\label{ODD}
 m(z)=-\sum\limits_{k=0}^\infty
\frac{\beta_k}{z^{2k+1}},\; \beta_0=1.
\end{equation}
Let $J$ be a bounded self-adjoint Jacobi matrix whose Weyl's
function coincides with $m(z).$ Then due to \eqref{ODD} we obtain
$b_1=0.$ By \eqref{connect} we have
\[
m_+(z,1)=\frac{1}{a^2_1}\left(-z-\frac{1}{m(z)}\right).
\]
It follows that $m_+(-z,1)=-m_+(z,1)$. Now \eqref{expansion} yields
$b_2=0$. Hence by induction and again using \eqref{connect} and
\eqref{expansion}, we obtain that $b_k=0$ for all $k.$
\end{proof}
\subsection{Mixed given data and uniqueness for finite self-adjoint Jacobi matrixces}
\label{subsec35}
 Let $J$ be an $n\times n$ self-adjoint Jacobi
matrix.  Following \cite{GS} we will consider $b$'s and $a$'s as a
single sequence $\{c_k\}_{k=1}^{2n-1}$, where $c_{2k-1}=b_k$ and
$c_{2k}=a_k$. The next theorem is established by H.~Hochstadt in
\cite{H3}.
\begin{theorem}
\label{HH} \cite{H3}. Let $J$ be an $n\times n$ self-adjoint Jacobi
matrix. Suppose that $c_{n+1},\ldots, c_{2n-1}$ are known as well as
the eigenvalues $z_1,\ldots,z_n$ of $J$. Then $c_1, \ldots, c_n$ are
uniquely determined.
\end{theorem}
F. ~Gesztesy and B.~Simon in \cite{GS} proved the following
generalization of the Hochstadt theorem.
\begin{theorem}
\label{GST} \cite{GS}. Let $J$ be an $n\times n$ self-adjoint Jacobi
matrix. Suppose that $c_{j+1},\ldots, c_{2n-1}$ are known as well as
$j$ of the eigenvalues. Then $c_1,\ldots, c_j$ are uniquely
determined.
\end{theorem}
Note that it is not necessary to know which of the $j$ eigenvalues one has,
and there may be no matrix consistent with the data as in Theorems
\ref{HH} and  \ref{GST} (see \cite{DN}).

Since later on we will consider a non-self-adjoint analog of the
Hochstadt and Gesztesy-Simon theorems, let us mention the key moment
of the proof of Theorem \ref{GST} in \cite{GS}.
 Apart from $m_+(z,k)$-functions the following $m_-$ - functions are
used in \cite{GS}:
\begin{equation}
\label{m-}
m_-(z,k)=((J_{[1,k-1]}-zI)^{-1}\delta_{k-1},\delta_{k-1}),\;
k=2,3,\ldots,n.
\end{equation}
It is established in \cite{GS} that \textit{for any eigenvalue
$\lambda_j$ of a self-adjoint $n\times n$ Jacobi matrix $J$ and for
any $k=1,\ldots,n$ the equality}
\begin{equation}
\label{mpm} m_-(\lambda_j,
k+1)=\left[a^2_km_+(\lambda_j,n)\right]^{-1}
\end{equation}
\textit{holds, where the equality in \eqref{mpm} includes the case
that both sides equal infinity.}
\section{Dissipative Jacobi matrices with a
rank one imaginary part} \label{sec4}
 \subsection{Finite and semi-infinite
dissipative Jacobi matrices}
 In the following we will consider the
Jacobi matrix of the form \eqref{22} and \eqref{222} but now we
suppose that condition \eqref{333} is fulfilled.
\begin{proposition} \label{LCHAR}
Let $J$ be a non-self-adjoint $n\times n$ or semi-infinite Jacobi
matrix of the form \eqref{22} with conditions \eqref{333} and
\eqref{Bound}. Then $J$ defines a prime dissipative operator with a
rank one imaginary part in $\dC^n$ $(l_2(\dN))$ and its Livsic
characteristic function takes the form
\begin{equation}
\label{LCHAR1} W(z)=1-2i\,\IM
\,b_1\left((J-zI)^{-1}\delta_1,\delta_1\right).
\end{equation}
\end{proposition}
\begin{proof}Let $J^*$ be the adjoint matrix to $J$ and let
\[
 J_R=\frac{1}{2}(J+J^*),\; J_I=\frac{1}{2i}(J-J^*)
\]
be Hermitian components of $J$. One has for $n\times n$ case
\begin{equation}
\label{5} J_R=\begin{pmatrix} \RE b_1 & a_1 & 0  &0  & 0 & \cdot &
\cdot & \cdot \\
a_1 & b_2 & a_2 & 0 &0& \cdot &
\cdot & \cdot \\
0    & a_2 & b_3 & a_3 &0& \cdot &
\cdot & \cdot \\
\cdot & \cdot & \cdot & \cdot &
\cdot & \cdot & \cdot &\cdot \\
\cdot & \cdot & \cdot & \cdot &
\cdot & \cdot & \cdot &a_{n-1} \\
\cdot & \cdot & \cdot & \cdot &\cdot & 0 & a_{n-1} & b_n
\end{pmatrix},\;
J_I=\begin{pmatrix}\IM b_1 & 0 & 0 &\cdot & \cdot &
\cdot & 0 \\
0 & 0 & 0& \cdot &
\cdot & \cdot& 0 \\
0    & 0 &0& \cdot &
\cdot & \cdot&0 \\
\cdot & \cdot & \cdot & \cdot &
\cdot & \cdot & \cdot \\
\cdot & \cdot & \cdot & \cdot &
\cdot & \cdot & \cdot \\
0 &0 &0&\cdot & \cdot&\cdot & 0
\end{pmatrix}
\end{equation}
and for semi-infinite case
\begin{equation}
\label{SEMHERM} J_R=\begin{pmatrix} \RE b_1 & a_1 & 0 &0   & 0 &
\cdot &
\cdot & \cdot \\
a_1 & b_2 & a_2 & 0 &0& \cdot &
\cdot & \cdot \\
0    & a_2 & b_3 & a_3 &0& \cdot &
\cdot & \cdot \\
\cdot & \cdot & \cdot & \cdot & \cdot & \cdot & \cdot &\cdot
\end{pmatrix},\;
J_I=\begin{pmatrix} \IM b_1 &  0 &0   & 0 & 0&\cdot &
\cdot & \cdot \\
0 & 0 & 0 & 0 &0&\cdot &
\cdot & \cdot \\
0    & 0 & 0 & 0 &0& \cdot &
\cdot & \cdot \\
\cdot & \cdot & \cdot & \cdot & \cdot & \cdot & \cdot &\cdot
\end{pmatrix},
\end{equation}
 and for every $ x\in \dC^n$ $(l_2(\dN))$
\[
2J_I x=\left( x, g\right)\, g,
\]
where
\begin{equation}
\label{canal}
  g=\sqrt{2\,\IM\, b_1}\,\delta_1.
\end{equation}
 By \eqref{Char} one obtains \eqref{LCHAR1}. From \eqref{coord} it follows that $J$ is prime.
\end{proof}
\begin{remark}
\label{RRM} One can easily show that a non-self-adjoint finite or
semi-infinite bounded tri-diagonal matrix $\cA$ with diagonal
entries $\{d_k\}$, sub-diagonal entries $\{c_k\}$, and
super-diagonal entries $\{\overline{c_k}$\} has a rank one imaginary
part if and only if there exists a diagonal nonzero entry $d_m$ such
that
\[
d_md_{m+1}=|c_m|^2,
\]
and remaining entries equal zero. It follows that the Jacobi matrix
$J$ with the conditions $\IM b_1=\cdots=\IM b_{m-1}=0$, $\IM b_m>
0$, $\IM b_{m+1}=\cdots=0$, $a_k>0$ for all $k$ is also dissipative
with a rank one imaginary part. But in this case such a matrix might
not be prime. For example the $3\times 3$ matrix of the form
\[
J=\begin{pmatrix}1&a_1&0\cr a_1&-1+it& a_2\cr 0 &a_2&1\end{pmatrix}
\]
is dissipative for $t>0$ and
\[
\ran J_I=\{\lambda \delta_2,\; \lambda\in\dC\}.
\]
But $\span\{\delta_2, J\delta_2, J^2\delta_2\}\ne \dC^3.$ So, this
matrix is not prime.
\end{remark}
 The next proposition easily follows from \eqref{kernel},
\eqref{SSSS}, and \eqref{SpFunc}.
\begin{proposition}
\label{Weyl} Let $J$ be a bounded non-self-adjoint finite  or
semi-infinite Jacobi matrix of the form \eqref{22} with condition
\eqref{333}. Then for the Weyl functions $m_J(z)$ and $m_{J_R}(z)$
of $J$ and $J_R$ and for the characteristic function $W(z)$ there
are the following relations:
\[
\begin{split}
&m_J(z)=\frac{i}{2\,\IM b_1}\,\left(W(z)-1\right),\\
&m_{J_R}(z)=\frac{i}{\IM b_1}\,\frac{W(z)-1}{W(z)+1},\\
&m_J(z)=\frac{\,im_{J_R}(z)}{i-\IM b_1\, m_{J_R}(z)},\\
&m_{J_R}(z)=\frac{\,im_{J}(z)}{i+\IM b_1\, m_{J}(z)}.\\
\end{split}
\]
Moreover, the function
\[
\cK_J(z,\xi):=\frac{m_J(z)-\ovl{m_J(\xi)}+2i\,\IM b_1
m_J(z)\ovl{m_J(\xi)}}{z-\bar\xi}
\]
is a nonnegative kernel, the function
\[
-\frac{1}{m_J(z)}+b_1
\]
is a Herglotz-Nevanlinna, and at infinity
\[
\frac{1}{a^2_1}\left(-\frac{1}{m_J(z)}+b_1-z\right)
\sim-\frac{1}{z}-\frac{b_2}{z^2}-\frac{b^2_2+a^2_2}{z^3}+O\left(\frac{1}{z^4}\right).
\]
The points of a non-real spectrum of $J$ are solutions of the
equation
\[
m_{J_R}(z)=\frac{i}{\IM b_1},\; \IM z>0.
\]
\end{proposition}
From  Proposition \ref{LCHAR}, Theorem \ref{model} and formulas
\eqref{V}, \eqref{FLT}, \eqref{wfuncAr}, \eqref{CONN1}, and
\eqref{SpFunc} we get the following statement.
\begin{proposition}
\label{spectrr} Let $J$ be a non-self-adjoint $n\times n$  Jacobi
matrix of the form \eqref{22} with conditions \eqref{333} and let
$z_1,z_2,\ldots,z_n$ be the eigenvalues of $J$ counting algebraic
multiplicity. Then the following formulas
\begin{equation}
\label{chweyl}
\begin{split}
&W(z)=\prod\limits_{k=1}^n\frac{z-\overline {z_k}}{z-{z_k}},\\
&m_J(z)=\frac{i}{2\,\IM b_1}\,
\frac{\prod\limits_{k=1}^n(z-\overline
{z_k})-\prod\limits_{k=1}^n(z-z_k)}{\prod\limits_{k=1}^n(z-{z_k})},\\
&m_{J_R}(z)=\frac{i}{\IM
b_1}\,\frac{\prod\limits_{k=1}^n(z-\overline
{z_k})-\prod\limits_{k=1}^n(z-z_k)}{\prod\limits_{k=1}^n(z-\overline
{z_k})+\prod\limits_{k=1}^n(z-{z_k})}
\end{split}
\end{equation}
hold.
\end{proposition}
\begin{proposition}
\label{dissi} Let $J$ be an $n\times n$  Jacobi matrix, $n\ge 3,$
with the conditions \eqref{333}. Then the matrices $J$ and
$J_{[1,n-2]}$ have no common eigenvalues and the relations
\begin{equation}
\label{invar} -\frac{a_1 a_2\cdots a_k}{(k-1)!}\,\IM
P^{(k-1)}_{k+1}(0)=\IM b_1
\end{equation}
hold for all $k=1,\ldots, n$, where the polynomials
$\{P_k(z)\}_{k=1}^{n+1}$ are defined by \eqref{SYSTEM}.
\end{proposition}
\begin{proof} Because all $J_{[1,k]}$ are dissipative and prime, their eigenvalues
 belong to $\dC_+$.
 Suppose that $z_0$ is an eigenvalue of $J$ and
$J_{[1,n-2]}$. Then by Proposition \ref{eigcores}
$P_{n-1}(z_0)=P_{n+1}(z_0)=0.$ From \eqref{SYSTEM} we have
\[
a_{n-1}P_{n-1}(z_0)+(b_n-z_0)P_n(z_0)=0.
\]
Since $\IM z_0>0$ and $b_n$ is real, we get that $P_{n}(z_0)=0$,
i.e. $z_0$ is a common eigenvalue for the matrices $J_{[1,n-2]}$ and
$J_{[1, n-1]}$. By Proposition \ref{COMMON} this is impossible.

Let $\{z_j\}_{j=1}^k$ be the roots of the polynomial $P_{k+1}(z)$.
Then $\IM z_j>0$. From \eqref{CHARPOL} it follows that
$$
\sum\limits_{j=1}^n z_j={\rm Sp} J_{[1,k]}.
$$
But $\IM {\rm Sp} J_{[1,k]}=\IM b_1.$ On the other hand
$$
\sum\limits_{j=1}^n z_j=-\frac{a_1 a_2\cdots a_k}{(k-1)!}\,
P^{(k-1)}_{k+1}(0).
$$
Therefore, \eqref{invar} holds.
\end{proof}

 If $J$ is a semi-infinite bounded Jacobi matrix of the form
\eqref{222} and with conditions \eqref{333}, then by Proposition
\ref{LCHAR} it determines a prime dissipative operator with rank one
imaginary part. Therefore its non-real spectrum is a denumerable set
of eigenvalues of finite algebraic multiplicities and from
\eqref{EIGENV}
\[
\sigma_p(J)=\left\{z\in\dC:
\sum\limits_{k=1}^\infty\left|P_k(z)\right|^2<\infty\right\},
\]
where $\{P_k(z)\}$ are corresponding polynomials defined by
\eqref{SYSTEM}. By Theorem \ref{Livsic1} the closure of the linear
span of all corresponding root subspaces coincides with $l_2(\dN)$
if and only if
\[
\sum\limits_{z\in\sigma_p(J)} \IM z=\IM b_1
\]
\begin{example}
\label{CHEB1}
 Let $J_l$ be the semi-infinite Jacobi matrix of the form
\begin{equation}
\label{jl} J_l=\begin{pmatrix} il & 1/2 & 0  &0  & 0 & \cdot &
\cdot & \cdot \\
1/2 &0 & 1/2 & 0 &0& \cdot &
\cdot & \cdot \\
0    & 1/2 & 0 & 1/2 &0& \cdot &
\cdot & \cdot \\
\cdot & \cdot & \cdot & \cdot &
\cdot & \cdot & \cdot &\cdot \\
\end{pmatrix},
\end{equation}
where $l>0$. The real part is the Jacobi matrix
\begin{equation}
\label{ho}
 H_0=\begin{pmatrix} 0 & 1/2 & 0  &0  & 0 & \cdot &
\cdot & \cdot \\
1/2 &0 & 1/2 & 0 &0& \cdot &
\cdot & \cdot \\
0    & 1/2 & 0 & 1/2 &0& \cdot &
\cdot & \cdot \\
\cdot & \cdot & \cdot & \cdot &
\cdot & \cdot & \cdot &\cdot \\
\end{pmatrix}.
\end{equation}
It is known \cite{Ber} that the Weyl function of $H_0$  is the
function
\[
m_{H_0}(z)=2(\sqrt{z^2-1}-z), \; z\in\dC\setminus [-1,1].
\]
The spectrum of $H_0$ is continuous and coincides with the interval
$[-1,1]$. Therefore, $||H_0||=1$. The Weyl function of $J_l$ is the
function
\[
m_{J_l}(z)=\frac{2i(\sqrt{z^2-1}-z)}{i-2l(\sqrt{z^2-1}-z)}.
\]
The non-real spectrum of $J_l$ we find solving the equation
\[
\sqrt{z^2-1}-z=\frac{i}{2l}.
\]
This equation has a unique solution from the open upper half-plane
\[
z=i\,\frac{4l^2-1}{4}
\]
for $l>1/2$ and has no solutions for $l\in (0,1/2]$. Thus if
$l\in(0,1/2]$ the spectrum of $J$ coincides with $[-1,1]$ and if
$l>1/2$ the spectrum of $J_l$ is $\{i(4l^2-1)/4\}\cup [-1,1].$

Note that
\[
\lim\limits_{x \uparrow -1}\frac{1}{m_{H_0}(z)}=-\lim\limits_{x
\downarrow 1}\frac{1}{m_{H_0}(z)}=\frac{1}{2}.
\]
Results from \cite{Kr}, \cite{ArTs2}, \cite{ArHdeS} yield that the
Jacobi matrices of the form
\[
 H_x=\begin{pmatrix} x & 1/2 & 0  &0  & 0 & \cdot &
\cdot & \cdot \\
1/2 &0 & 1/2 & 0 &0& \cdot &
\cdot & \cdot \\
0    & 1/2 & 0 & 1/2 &0& \cdot &
\cdot & \cdot \\
\cdot & \cdot & \cdot & \cdot &
\cdot & \cdot & \cdot &\cdot \\
\end{pmatrix}
\]
determine contractive Jacobi operators in $l_2(\dN)$ if and only if
$|x|\le 1/2$.
\end{example}
\begin{example}
\label{Cheb2} Let $\wh{J}_l$ be the semi-infinite Jacobi matrix of
the form
\begin{equation}
\label{hatjl} \wh J_l=\begin{pmatrix} il & 1/\sqrt{2} & 0  &0  & 0 &
\cdot &
\cdot & \cdot \\
1/\sqrt{2} &0 & 1/2 & 0 &0& \cdot &
\cdot & \cdot \\
0    & 1/2 & 0 & 1/2 &0& \cdot &
\cdot & \cdot \\
\cdot & \cdot & \cdot & \cdot &
\cdot & \cdot & \cdot &\cdot \\
\end{pmatrix},
\end{equation}
where $l>0$. The real part of $\wh{J}_l$ is the Jacobi matrix
\begin{equation}
\label{hatho}
 \wh{H}_0=\begin{pmatrix} 0 & 1/\sqrt{2} & 0  &0  & 0 & \cdot &
\cdot & \cdot \\
1/\sqrt{2} &0 & 1/2 & 0 &0& \cdot &
\cdot & \cdot \\
0    & 1/2 & 0 & 1/2 &0& \cdot &
\cdot & \cdot \\
\cdot & \cdot & \cdot & \cdot &
\cdot & \cdot & \cdot &\cdot \\
\end{pmatrix}.
\end{equation}
It is known \cite{Ber} that the Weyl function of $\wh{H}_0$  is the
function
\[
m_{\wh{H}_0}(z)=-\frac{1}{\sqrt{z^2-1}}, \; z\in\dC\setminus [-1,1].
\]
The spectrum of $\wh H_0$ is continuous and coincides with the
interval $[-1,1]$. Therefore $||\wh H_0||=1.$ The Weyl function of
$\wh J_l$ is the function
\[
m_{\wh J_l}(z)=-\frac{i}{i\sqrt{z^2-1}+l}.
\]
The non-real spectrum of $\wh J_l$ we find solving the equation
\[
\sqrt{z^2-1} =il.
\]
This equation has a unique solution from the open upper half-plane
\[
z=i\sqrt{l^2-1}
\]
for $l>1$ and has no solutions for $l\in (0,1]$. Thus if $l\in(0,1]$
the spectrum of $\wh J_l$ coincides with $[-1,1]$ and if $l>1$ the
spectrum of $J_l$ is $\{i\sqrt{l^2-1}\}\cup [-1,1].$

Note that
 \[
 \lim\limits_{x \uparrow -1}{|m_{\wh
H_0}(x)|}=\lim\limits_{x \downarrow 1}{ |m_{\wh H_0}(x)|}=+\infty.
\]
 From \cite{Kr} it follows that the Jacobi matrix $\wh H_0$ is a
unique contractive matrix in $l_2(\dN)$ among Jacobi matrices of the
form
\[
 \wh{H}_x=\begin{pmatrix} x & 1/\sqrt{2} & 0  &0  & 0 & \cdot &
\cdot & \cdot \\
1/\sqrt{2} &0 & 1/2 & 0 &0& \cdot &
\cdot & \cdot \\
0    & 1/2 & 0 & 1/2 &0& \cdot &
\cdot & \cdot \\
\cdot & \cdot & \cdot & \cdot &
\cdot & \cdot & \cdot &\cdot \\
\end{pmatrix}.
\]
\end{example}

\subsection{The $m$-functions of a finite dissipative Jacobi matrix with a rank
one imaginary part} \label{subsec4}
 The next proposition establishes
the analogues of the relations \eqref{mpm} for a dissipative Jacobi
matrix.
\begin{theorem}
\label{Green} Let $J$ be an $n\times n$ Jacobi matrix with
conditions \eqref{333}. Let $m_J(z)$ be the Weyl function of $J$
defined by \eqref{SpFunc} and let the functions $m_+(z,j)$ and
$m_-(z,j)$ be defined by \eqref{mk} and \eqref{m-}. If $z_0$ is the
eigenvalue of $J$ of the algebraic multiplicity $l$ then  for all
$j=1,\ldots ,n$ the following relations hold
\begin{equation}
\label{new}
\left(\frac{1}{m_-(z,j+1)}\right)^{(p)}\Biggr|_{z=z_0}=a^2_j\,m^{(p)}_+(z_0,j),\;
p=0,\ldots,l-1,\;
\end{equation}
\end{theorem}
\begin{proof}
Represent the matrix $J$ in the block-matrix form
\[
J=\begin{pmatrix}J_{[1,j]}&A_{12}\cr
A_{21}&J_{[j+1,n]}\end{pmatrix},
\]
where $A_{12}$ is $j\times (n-j)$ matrix  of the form
\[
A_{12}=\begin{pmatrix} 0&0&\cdots&0\cr 0&0&\cdots&0\cr
\cdots&\cdots&\cdots&\cdots\cr0&0&\cdots&0\cr a_j&0&\cdots&0
\end{pmatrix}
\]
and $A_{21}=A^*_{12}$ (is $(n-j)\times j$ matrix). Let
$\cH_1=\dC^j$. Then from \eqref{ShF} it follows that
\[
P_{\cH_1}(J-z)^{-1}\uphar\cH_1=\left(J_{[1,j]}-zI-A_{12}(J_{[j+1,n]}-zI)^{-1}A_{21}\right)^{-1}
\]
for $z\in\rho(J)\cap\rho\left(J_{[j+1,n]}\right)$. Clearly
\[
A_{12}(J_{[j+1,n]}-zI)^{-1}A_{21}=\begin{pmatrix} 0&0&\cdots&0&0\cr
0&0&\cdots&0&0\cr
\cdots&\cdots&\cdots&\cdots&\cdots\cr0&0&\cdots&0&0\cr 0&0&\cdots&0&
a^2_jm_+(z,j)
\end{pmatrix}.
\]
Now from \eqref{RESFOR} we get
\begin{equation}
\label{Green1}
\left((J-zI)^{-1}\delta_j,\delta_j\right)=\frac{m_-(z,j+1)}{1-a^2_j\,m_+(z,j)\,m_-(z,j+1)},
\;j=1,2,\dots, n.
\end{equation}
Let $j\ge 2.$ By Cramer's rule we have
\[
\left((J-zI)^{-1}\delta_j,\delta_j\right)=\frac{{\rm
det}(J_{[1,j-1]}-zI)\,{\rm
 det}(J_{[j+1,n]}-zI)}{{\rm
det}(J-zI)}
\]
and
\[
m_-(z,j+1)=\frac{{\rm det}(J_{[1,j-1]}-zI)}{{\rm det}(J_{[1,j]}-zI)}
\]
for $j=2,\dots, n$. From \eqref{Green1} it follows that
\[
\frac{{\rm
 det}(J_{[j+1,n]}-zI)}{{\rm
det}(J-zI)}=\frac{1}{{\rm det}(J_{[1,j]}-zI) -a^2_j\,m_+(z,j)\,{\rm
det}(J_{[1,j-1]}-zI)}
\]
Let $z_0$ be an eigenvalue of $J$ of the algebraic multiplicity $l$.
Because $J$ is dissipative, we have $\IM z_0>0$. Note that
$J_{[j+1,n]}$ is self-adjoint Jacobi matrix, the Herglotz-Nevanlinna
function $m_+(z,j)$ is holomorphic at $z_0$ and $\IM m_+(z_0,j)> 0$.
It follows that the function
\[
\frac{{\rm
 det}(J_{[j+1,n]}-zI)}{{\rm
det}(J-zI)}
\]
has the pole at $z_0$ of order $l$. Therefore
\[
\lim\limits_{z\to z_0}(z-z_0)^p\,\frac{{\rm
 det}(J_{[j+1,n]}-zI)}{{\rm
det}(J-zI)}=\infty,\; p=0,1,\ldots l-1
\]
and
\[
\lim\limits_{z\to z_0}(z-z_0)^l\,\,\frac{{\rm
 det}(J_{[j+1,n]}-zI)}{{\rm
det}(J-zI)} \quad\mbox{is finite}.
\]
This yields
\begin{equation}
 \label{limitsp}
\lim\limits_{z\to z_0} \frac{(z-z_0)^p}{{\rm det}(J_{[1,j]}-zI)
-a^2_j\,m_+(z,j)\,{\rm det}(J_{[1,j-1]}-zI)}=\infty
\end{equation}
for $p=0,1,\ldots,l-1$. In particular, for $p=0$ we obtain
\begin{equation}
\label{hol1} {\rm det}(J_{[1,j]}-z_0I) -a^2_j\,m_+(z_0,j)\,{\rm
det}(J_{[1,j-1]}-z_0I)=0.
\end{equation}
By Proposition \ref{COMMON} the matrices $J_{[1,j-1]}$ and $J_{[1,
j]}$ have no common eigenvalues and since $m_+(z_0,j)\ne 0$,
\eqref{hol1} implies that
\[
{\rm det}(J_{[1,j]}-z_0I)\ne 0,\; {\rm det}(J_{[1,j-1]}-z_0I)\ne 0
\]
and
\begin{equation}
\label{hol} \frac{1}{m_-(z_0,j+1)}=\frac{{\rm
det}(J_{[1,j]}-z_0I)}{{\rm
det}(J_{[1,j-1]}-z_0I)}=a^2_j\,m_+(z_0,j).
\end{equation}
In some neighborhood of $z_0$ we can rewrite \eqref{limitsp} as
\[
\lim\limits_{z\to z_0} \frac{1}{{\rm
det}(J_{[1,j-1]}-zI)}\,\frac{(z-z_0)^p}{(m_-(z,j+1))^{-1}
-a^2_j\,m_+(z,j)}=\infty
\]
for $p=0,1,\ldots, l-1$. It follows that the number $z_0$ is a zero
of the function
\[
\frac{1}{m_-(z,j+1)}-a^2_j\,m_+(z,j)
\]
of order $l$. Thus, we get \eqref{new} for $j\ge 2$.

Let $j=1$. Since $m_-(z,2)=(b_1-z)^{-1}$, relation \eqref{Green1}
takes the form
\[
m_J(z)=\frac{1}{b_1-z -a^2_1 m_+(z,1)}.
\]
If $z_0$ is the eigenvalue of $J$ of the algebraic multiplicity $l$,
then due to \eqref{chweyl} the Weyl function $m_J(z)$ has  a pole of
order $l$ at $z_0$. In the same manner as the above we get that $z_0$ is a zero of
order $l$ for the function
\[
b_1-z-a^2_1 m_+(z,1)=\frac{1}{m_-(z,2)}-a^2_1 m_+(z,1).
\]
Therefore, \eqref{new} holds and for $j=1$.
\end{proof}
 Since the
integer $j$ is an arbitrary, we obtain as a by-product the following
corollary.
\begin{corollary}
\label{COROL} Let $J$ be an $n\times n$ Jacobi matrix with the
conditions \eqref{333}. Then the matrices $J_{[1,j]}$ and $J$ have
no common eigenvalues for every $j=1,2,\ldots, n-1$.
\end{corollary}

\section{Inverse spectral problems for finite Jacobi matrices with a
rank one imaginary part} \label{sec5}
\subsection{Reconstruction of a finite
dissipative Jacobi matrix from its eigenvalues}
 Let
$J$ be a non-self-adjoint $n\times n$ Jacobi matrix of the form
\eqref{22} satisfying conditions \eqref{333}. Then by Proposition
\ref{LCHAR} the corresponding operator in $\dC^n$ is a prime
dissipative operator with a rank one imaginary part. Therefore the
matrix $J$ has only non-real eigenvalues with positive imaginary
parts. The next theorem establishes that an arbitrary $n$ non-real
numbers counting multiplicity from the open upper half-plane
determine uniquely some dissipative $n\times n$ Jacobi matrix with a
rank one imaginary part.
\begin{theorem}
 \label{TT1}
 Suppose that $z_1, \ldots, z_n$ are not
necessarily distinct complex numbers with positive imaginary parts.
Then there exists a unique $n\times n $ Jacobi matrix with entries
satisfying conditions \eqref{333} whose eigenvalues counting
algebraic multiplicity coincide with $\{z_k\}_{k=1}^n.$
\end{theorem}
\begin{proof} Let
\begin{equation}
\label{CCC}
 W(z)=\prod\limits_{k=1}^n\frac{z-{\bar z}_k}{z-z_k}.
\end{equation}
Then
\[
\lim\limits_{z\to \infty}z(W(z)-1)=2i\sum_{k=1}^n\IM z_k.
\]
Let
\[
c=\sum_{k=1}^n\IM z_k.
\]
Define
\begin{equation}
\label{m+0}
 m_+(z,0)=\frac{i}{c}\,\,\frac{W(z)-1}{W(z)+1}.
\end{equation}
The Herglotz-Nevanlinna function $m_+(z,0)$ has the expansion at
infinity
\[
m_+(z,0)\sim-\frac{1}{z}-\frac{b}{z^2}-\frac{b^2+a^2_1}{z^3}+O\left(\frac{1}{z^4}\right),
\]
and determines a probability measure supported at some $n$ points.
By Theorem \ref{recov} there exists a unique self-adjoint $n\times n$
Jacobi matrix
\[
H=\begin{pmatrix} b & a_1 & 0    & 0 & \cdot &
\cdot & \cdot \\
a_1 & b_2 & a_2 & 0 & \cdot &
\cdot & \cdot \\
0    & a_2 & b_3 & a_3 & \cdot &
\cdot & \cdot \\
\cdot & \cdot & \cdot & \cdot &
\cdot & \cdot & \cdot \\
\cdot & \cdot & \cdot & \cdot &
\cdot & \cdot & a_{n-1} \\
\cdot & \cdot & \cdot & \cdot & 0 & a_{n-1} & b_n
\end{pmatrix}
\]
satisfying conditions \eqref{33} such that
\[
m_+(z,0)=\left((H-zI)^{-1}\delta_{1},\delta_{1}\right).
\]
Let
\[
J=\begin{pmatrix} b+ic & a_1 & 0    & 0 & \cdot &
\cdot & \cdot \\
a_1 & b_2 & a_2 & 0 & \cdot &
\cdot & \cdot \\
0    & a_2 & b_3 & a_3 & \cdot &
\cdot & \cdot \\
\cdot & \cdot & \cdot & \cdot &
\cdot & \cdot & \cdot \\
\cdot & \cdot & \cdot & \cdot &
\cdot & \cdot & a_{n-1} \\
\cdot & \cdot & \cdot & \cdot & 0 & a_{n-1} & b_n
\end{pmatrix}
=H+i\begin{pmatrix} c & 0 & 0    & 0 & \cdot &
\cdot & \cdot \\
0 & 0 & 0 & 0 & \cdot &
\cdot & \cdot \\
0    & 0 &0 & 0 & \cdot &
\cdot & \cdot \\
\cdot & \cdot & \cdot & \cdot &
\cdot & \cdot & \cdot \\
\cdot & \cdot & \cdot & \cdot &
\cdot & \cdot & \cdot \\
\cdot & \cdot & \cdot & \cdot & 0 & 0 & 0
\end{pmatrix}.
\]
It follows that
\[
2\,J_I\, x =\left( x, g\right), \; x\in\dC^n,
\]
where $ g=\sqrt{2c}\delta_1\in\dC^n$.

 By Proposition \ref{LCHAR} for the characteristic
function of $J$ we obtain
\[
w(z)=1-2i\,c\left((J-zI)^{-1}\delta_1,\delta_1\right).
\]
Its fractional-linear transformation
\[
v(z)=i\,\frac{w(z)-1}{w(z)+1}
\]
has the form
\[
v(z)=\frac{1}{2}\left((H-zI)^{-1}g,
g\right)=c\,\left((H-zI)^{-1}\delta_1,\delta_1\right)=c\,m_+(z,0).
\]
From \eqref{m+0} we get
\[
v(z)=i\,\frac{W(z)-1}{W(z)+1}.
\]
Therefore $w(z)=W(z)$. By Theorem \ref{model} the eigenvalues of $J$
counting algebraic multiplicity coincide with $\{z_k\}$.
\end{proof}
\begin{example}
Let us construct $3\times 3$ dissipative Jacobi matrix with
eigenvalues $z_1=i$ of multiplicity 2 and $z_2=2i$ (of multiplicity
1). By Theorem \ref{Livsic} the characteristic function of $J$ is
\[
W(z)=\left(\frac{z+i}{z-i}\right)^2\frac{z+2i}{z-2i}.
\]
Then
\[
\lim\limits_{z\to\infty}z(W(z)-1)=4i.
\]
Note  that $2\IM z_1+ \IM z_2=4.$ Let
 \[
\begin{split}
&m_+(z,0)=\frac{i}{4}\,\frac{W(z)-1}{W(z)+1}=\frac{i}{4}\,\frac{(z+i)^2(z+2i)-(z-i)^2(z-2i)}
{(z+i)^2(z+2i)+(z-i)^2(z-2i)}=\\
&\qquad=\frac{-2z^2+1}{2z^3-10 z}\,.
\end{split}
\]
Then
\[
m(z)=\frac{-1}{z-\displaystyle\frac{9/2}{z-\displaystyle\frac{\mathstrut1/2}{z}}}.
\]
It follows that $b_1=4i,$ $b_2=b_3=0$, $a_1=\sqrt{9/2}$,
$a_2=\sqrt{1/2}$. The Jacobi matrix takes the form
\[
J=\begin{pmatrix} 4i&\frac{3}{\sqrt{2}}&0\cr
\frac{3}{\sqrt{2}}&0&\frac{1}{\sqrt{2}}\cr 0 &\frac{1}{\sqrt{2}}&0
\end{pmatrix}.
\]
\end{example}
\begin{example}
In order to construct the dissipative $n\times n$ Jacobi matrix with
the eigenvalue $z_0=x_0+i y_0$ of algebraic multiplicity $n$
($y_0>0$) it is sufficient to construct the $ n\times n$ Jacobi
matrix $J_n$ with the eigenvalue $z_0=i$ of algebraic multiplicity
$n$. The characteristic function of such a matrix is
\[
W_n(z)=\left(\frac{z+i}{z-i}\right)^n.
\]
The corresponding Weyl function of the real part is
\[
m_n(z)=\frac{i}{n}\,\,\frac{(z+i)^n-(z-i)^n}{(z+i)^n+(z-i)^n}\,.
\]
\end{example}
Because ${\rm Sp} J=in$, we get $\IM b_1=n$. Expanding the function
$m_n(z)$ via continued fraction we obtain the real part of $J_n$. In
particular
\[
\begin{split}
& m_2(z)=\frac{-1}{z-\displaystyle\frac{1}{z}},\;
m_3(z)=\frac{-1}{z-\displaystyle\frac{8/3}{z-\displaystyle\frac{\mathstrut1/3}{z}}}\,,\\
&m_4(z)=\frac{-1}{z+\displaystyle\frac{-5}{z+\displaystyle\frac{\mathstrut
-4/5}{z+\displaystyle\frac{\mathstrut-1/5}{z}}}},\;
m_5(z)=\frac{-1}{z+\displaystyle\frac{-8}{z+\displaystyle\frac{\mathstrut
-7/5}{z+\displaystyle\frac{\mathstrut-16/35}{z+\displaystyle\frac{\mathstrut-1/7}{z}}}}}
\end{split}.
\]
Therefore
\[
J_2=\begin{pmatrix}2i& 1\cr 1&0 \end{pmatrix},\;
J_3=\begin{pmatrix}3i& \frac{2\sqrt{2}}{\sqrt{3}}&0\cr
\frac{2\sqrt{2}}{\sqrt{3}}&0&\frac{1}{\sqrt{3}}\cr
0&\frac{1}{\sqrt{3}}&0
\end{pmatrix},
\]
\[
J_4=\begin{pmatrix}4i& \sqrt{5}&0&0\cr
\sqrt{5}&0&\frac{2}{\sqrt{5}}&0\cr
0&\frac{2}{\sqrt{5}}&0&\frac{1}{\sqrt{5}}\cr
0&0&\frac{1}{\sqrt{5}}&0
\end{pmatrix},\;
J_5=\begin{pmatrix}5i& \frac{2\sqrt{2}}{\sqrt{5}}&0&0&0\cr
\frac{2\sqrt{2}}{\sqrt{5}}&0&\frac{\sqrt{7}}{\sqrt{5}}&0&0\cr
0&\frac{\sqrt{7}}{\sqrt{5}}&0&\frac{4}{\sqrt{35}}&0\cr
0&0&\frac{4}{\sqrt{35}}&0&\frac{1}{\sqrt{7}}&\cr
0&0&0&\frac{1}{\sqrt{7}}&0
\end{pmatrix}.
\]
Note that the eigenvalues of the real part of $J_n$ are solutions of
the equations
\[
\left(\frac{z+i}{z-i}\right)^n=-1,
\]
that are real numbers
$$\lambda_k=-\cot\left(\frac{\pi+2\pi k}{2n}\right),\; k=0,1\ldots, n-1.$$
\subsection{The connection between triangular and Jacobi forms of a prime
dissipative operator with a rank one imaginary part}
 Let $z_1,z_2,\ldots,z_n$ be not
necessarily distinct complex number with positive imaginary parts.
According to Theorem \ref{TT1} there exists a unique $n\times n $
Jacobi matrix
\[
J=\begin{pmatrix} b_1 & a_1 & 0    & 0 & \cdot &
\cdot & \cdot \\
a_1 & b_2 & a_2 & 0 & \cdot &
\cdot & \cdot \\
0    & a_2 & b_3 & a_3 & \cdot &
\cdot & \cdot \\
\cdot & \cdot & \cdot & \cdot &
\cdot & \cdot & \cdot \\
\cdot & \cdot & \cdot & \cdot &
\cdot & \cdot & a_{n-1} \\
\cdot & \cdot & \cdot & \cdot & 0 & a_{n-1} & b_n
\end{pmatrix}
\]
 with entries satisfying conditions \eqref{333}
whose eigenvalues coincide with $\{z_k\}_{k=1}^n$ counting algebraic
multiplicity. On the other hand such a matrix by Theorem \ref{model}
is unitary equivalent to the triangular matrix of the form
\begin{equation}
\label{10} \vec{J}=\begin{pmatrix}z_1 & i\beta_1\beta_2 & \cdot &
\cdot & \cdot & i\beta_1\beta_n \\
0 & z_2 &  \cdot &
\cdot & \cdot & i\beta_2\beta_n \\
\cdot & \cdot & \cdot & \cdot &
\cdot & \cdot  \\
0 & 0 & \cdot & \cdot &
\cdot & z_n \\
\end{pmatrix},
\end{equation}
where
\begin{equation}
\label{9} z_k=\alpha_k+ i\frac{\beta^2_k}{2}=\IM z_k, \; \beta_k>0,
\; k=1,\ldots,n.
\end{equation}
From \eqref{10} it follows that
\[
2\,\IM\vec{J}=\begin{pmatrix} \beta^2_1 & \beta_1\beta_2 & \cdot &
\cdot & \cdot & \beta_1\beta_n \\
\beta_1\beta_2 & \beta^2_2 &  \cdot &
\cdot & \cdot & \beta_2\beta_n \\
\cdot & \cdot & \cdot & \cdot &
\cdot & \cdot  \\
\beta_1\beta_n & \beta_2\beta_n & \cdot & \cdot &
\cdot & \beta^2_n \\
\end{pmatrix}
\]
and
\begin{equation}
\label{12} 2\,\IM\vec J=\left(\cdot, \vec g\right)\, \vec g,
\end{equation}
where
\begin{equation}
\label{13}
 \vec g=\begin{pmatrix}\beta_1\\
\beta_2\\
\cdot\\
\cdot\\
\cdot \\
\beta_n
\end{pmatrix}.
\end{equation}
Let $U$ be unitary matrix
\[
U=\begin{pmatrix} u_{11} & u_{12} & \cdot &
\cdot & \cdot & u_{1n} \\
u_{21} & u_{22} &  \cdot &
\cdot & \cdot &u_{2n} \\
\cdot & \cdot & \cdot & \cdot &
\cdot & \cdot  \\
u_{n1} & u_{n2}& \cdot & \cdot &
\cdot & u_{nn} \\
\end{pmatrix}
\]
such that
\begin{equation}
\label{15} \left\{\begin{array}{l} U\,J=\vec J\,U\\
U{ g}=\vec{ g}
\end{array}\right.,
\end{equation}
where $g=\sqrt{2\IM b_1}\,\delta_1.$ Next we present \textit{an
algorithm which allows to find the Jacobi matrix $J$ and the unitary
matrix $U$ satisfying \eqref{15}}. From \eqref{15} it follows
\begin{equation}
\label{16}
 U\,J^k{ g}=(\vec{J})^k\vec{ g},\; k=1,2,...,n.
\end{equation}
 Since $U$ is the unitary matrix, we have $ \|{ g}\|^2=\|\vec{ g}\|^2
 $,
and taking into account \eqref{canal} and \eqref{13} we get
\[
\|{ g}\|^2=2\IM\,
b_1=\beta^2_1=\beta^2_2+...+\beta^2_n=\|\vec{g}\|^2.
\]
Thus
\begin{equation}
\label{18}  \IM\, b_1=\sum\limits_{k=1}^n\IM\,z_k.
\end{equation}
 From $U{ g}=\vec{ g}$ it follows
\begin{equation}
\label{20} u_{k1}=\frac{\sqrt{\IM\,
z_k}}{\sqrt{\sum\limits_{j=1}^n\IM\,z_j}},\;\, k=1,2,..., n.
\end{equation}
Relations \eqref{15} and \eqref{16} yield
\[
\left(\vec J\vec{ g},\vec{ g}\right)=\left(UJ{ g},U{
g}\right)=\left(J{ g},{ g}\right).
\]
Taking into account \eqref{1}, \eqref{canal}, \eqref{10} and
\eqref{13}, we get
\begin{equation}
\label{n22}
J{ g}=\begin{pmatrix}b_1\sqrt{2\,\IM\, b_1}\\
a_1\sqrt{2\,\IM\, b_1}\\
0\\
\cdot\\
\cdot\\
\cdot \\
0
\end{pmatrix},\; \vec J\vec{ g}=\begin{pmatrix}\left(z_1+i\sum\limits_{k=2}^n\beta^2_k\right)\beta_1\\
\left(z_2+i\sum\limits_{k=3}^n\beta^2_k\right)\beta_2\\
\cdot\\
\cdot\\
\cdot \\
z_n\beta_n
\end{pmatrix}.
\end{equation}
Therefore
\[
\begin{split}
&\left(J{ g},{ g}\right)=2\, b_1\,\IM b_1,\\
&\left(\vec J\vec{ g},\vec{
g}\right)=\sum\limits_{j=1}^n\left(z_j+i\sum\limits_{k=j+1}^n\beta^2_k\right)\beta^2_j
.
\end{split}
\]
 Since
$\left(J{g},{ g}\right)=\left(\vec J\vec{ g},\vec{ g}\right),$ from
\eqref{9} and \eqref{18}we get
\[
b_1=\frac{\sum\limits_{j=1}^n\left(z_j+i\sum\limits_{k=j+1}^n\beta^2_k\right)\IM
z_j} {\sum\limits_{k=1}^n\IM z_k},\; \RE
b_1=\frac{\sum\limits_{k=1}^n\RE z_k\,\IM z_k}
{\sum\limits_{k=1}^n\IM z_k}.
\]
Because $U$ is unitary and $UJ{ g}=\vec J\vec{ g}$, we get
$\left\|J{ g}\right\|^2=\left\|\vec J\vec{ g}\right\|^2.
$ This equality and \eqref{n22} yield
\[
|b_1|^2\,(2\,\IM b_1)+a^2_1(2\,\IM
b_1)=\sum\limits_{j=1}^n\left|z_j+i\sum\limits_{k=j+1}^n\beta^2_k\right|^2\beta^2_j
\]
and hence
\[
a_1=\left(\frac{\sum\limits_{j=1}^n\left|z_j+i\sum\limits_{k=j+1}^n\beta^2_k\right|^2\IM
z_j}{\sum\limits_{k=1}^n\IM z_k}\,-
\,\frac{\left|\sum\limits_{j=1}^n\left(z_j+i\sum\limits_{k=j+1}^n\beta^2_k\right)\IM
z_j\right|^2}{\left(\sum\limits_{k=1}^n\IM
z_k\right)^2}\right)^{1/2}.
\]

Recall that by \eqref{coord}
$$(J^{m-1}{g})_m=a_{m-1}a_{m-2}\cdot\cdot\cdot a_1\sqrt{2\,\IM
b_1}$$
and
 $$(J^{m-1}{g})_{m+1}=(J^{m-1}{
g})_{m+2}=\cdots=(J^{m-1}{g})_n=0$$ for $m=1,\ldots, n$. Let $m\ge
2$. Suppose that $(J^{m-1}{ g})_k$ are already known, where
$k=1,2..., m-1$. Then from \eqref{1} we obtain
\[
J^m{g}=\begin{pmatrix}b_1(J^{m-1}{g})_1+a_1(J^{m-1}{ g})_2\\
a_1 (J^{m-1}{g})_1+b_2(J^{m-1}{g})_2+a_2(J^{m-1}{g})_3\\
a_2 (J^{m-1}{g})_2+b_3(J^{m-1}{g})_3+a_3(J^{m-1}{g})_4
\\
\cdot\\
\cdot\\
\cdot \\
a_{m-1}(J^{m-1}{g})_{m-1}+b_m(J^{m-1}{g})_m\\
a_m(J^{m-1}{g})_m \\0\\
\cdot\\
\cdot\\
\cdot\\
0
\end{pmatrix}.
\]
It follows that
\[
\begin{split}
&\left(J^m{ g}, J^{m-1}{g}\right)_{\dC^n}=\sum\limits_{j=1}^n (J^m{
g})_j\,\ovl{(J^{m-1}{ g})_j}=\sum\limits_{j=1}^{m-1}
(J^m{g})_j\,\ovl{(J^{m-1}{
g})_j}+\\
&+\left(a_{m-1}(J^{m-1}{ g})_{m-1}+b_m(J^{m-1}{
g})_m\right)\ovl{(J^{m-1}{ g})_m}.
\end{split}
\]
From
\[
\left(J^m{ g}, J^{m-1}{ g}\right)=\left((\vec J)^m\vec{ g}, (\vec
J)^{m-1}\vec{ g}\right)
\]
and $(J^{m-1}{ g})_m=a_{m-1}a_{m-2}\cdot\cdot\cdot a_1\sqrt{2\,\IM
b_1}$ we get the linear equation with respect to $b_m$:
\[
\begin{split}
&\sum\limits_{j=1}^{m-1} (J^m{g})_j\,\ovl{(J^{m-1}{g})_j}
+\left(a_{m-1}(J^{m-1}{g})_{m-1}+b_m(J^{m-1}{
g})_m\right)\ovl{(J^{m-1}{ g})_m}=\\
&\qquad\qquad\qquad=\left((\vec J)^m\vec{ g}, (\vec J)^{m-1}\vec{
g}\right).
\end{split}
\]
This equation can be solved since the coefficient $(J^{m-1}{ g})_m$
is nonzero by \eqref{coord}. In order to find $a_m$ we use the
relation
\[
\left\|J^m{g}\right\|^2=\left\|(\vec J)^m\vec{ g} \right\|^2
\]
which takes the form
\[
\begin{split}
&\sum\limits_{j=1}^{m-1}\left|(J^m{
g})_j\right|^2+\left|a_{m-1}(J^{m-1}{g})_{m-1}+b_m(J^{m-1}{
g})_m\right|^2+\\
&\qquad\qquad\qquad+a^2_m[J^{m-1}{\bf g}]^2_m=\left\|(\vec J)^m\vec{
g} \right\|^2.
\end{split}
\]
Solving for $a^2_m$ we can find $a_m$.

According to \eqref{20} the elements $\{u_{k1}\}_{k=1}^n$ are
expressed by means of $z_1,z_2,\ldots,z_n$. The equality $UJ{
g}=\vec J\vec{ g}$ gives the following linear system with respect to
the second column $\{u_{k2}\}_{k=1}^n$ of the matrix $U$:
\[
\left\{
\begin{split}
&u_{11}(J{ g})_1+u_{12}( J{g})_2=(\vec J\vec{ g})_1\\
&u_{21}(J{ g})_1+u_{22}(J{ g})_2=(\vec
J\vec{ g})_2\\
&\ldots\quad\ldots\\
&u_{n1}(J{ g})_1+u_{n2}(J{ g})_2=(\vec
J\vec{ g})_n\\
\end{split}
\right..
\]
Since $(J{ g})_2=a_1\sqrt{2\,\IM b_1}\ne 0$, one can find
$\{u_{k2}\}_{k=1}^n.$ By induction the equality $UJ^{m-1}{ g}=(\vec
J)^{m-1}\vec{ g}$ enables us to find $\{u_{k\,m}\}_{k=1}^n$ for
$m\le n.$

Finally the relation $J=U^{-1}\vec JU$ and already known entries
$\{a_k\}_{k=1}^{n-1}$ and $\{b_k\}_{k=1}^{n-1}$ allow us to find the
entry $b_n.$

In conclusion we give formulas for the reconstruction of $2\times 2$
dissipative Jacobi matrix $J$ with a rank one imaginary part from
its eigenvalues $z_1 , z_2$, and the corresponding unitary matrix
$U$:
\[
J=\begin{pmatrix} \frac{\RE z_1\IM z_1+\RE z_2 \IM z_2}{\IM z_1+\IM
z_2}+i(\IM z_1+\IM z_2)&\sqrt{\left(\left(\frac{\RE z_1-\RE z_2}{\IM
z_1+\IM z_2}\right)^2+1\right)\,\IM z_1\,\IM z_2}\cr
\sqrt{\left(\left(\frac{\RE z_1-\RE z_2}{\IM z_1+\IM
z_2}\right)^2+1\right)\,\IM z_1\,\IM z_2} & \frac{\RE z_1\,\IM
z_2+\RE z_2\,\IM z_1}{\IM z_1+\IM z_2}
\end{pmatrix},
\]
\[
U=\begin{pmatrix}\sqrt{\displaystyle\frac{\IM z_1}{\IM z_1+\IM
z_2}}&\sqrt{\displaystyle\frac{\IM z_2}{\IM z_1+\IM z_2}}\cr
 \sqrt{\displaystyle\frac{\IM z_2}{\IM z_1+\IM z_2}}&-\sqrt{\displaystyle\frac{\IM z_1}{\IM z_1+\IM
 z_2}}\end{pmatrix}.
\]
\begin{remark}The presented algorithm of reconstruction of the
unique dissipative non-self-adjoint Jacobi matrix with a rank one
imaginary part having given non-real numbers as its eigenvalues
allows to recover the set of tri-diagonal matrices with the same
non-real eigenvalues.
\end{remark}

\subsection{A non-self-adjoint analog of the Hochstadt
and Gesztesy--Simon theorems} \label{subsec5} The next theorem is a
non-self-adjoint analog the of Hochstadt \cite{H3} and Gesztesy
-Simon \cite{GS} results (see Theorems \ref{HH} and \ref{GST}).
\begin{theorem}
\label{AnHochs}Let $J$ be an $n\times n$ Jacobi matrix with
conditions \eqref{333}. Suppose that the eigenvalues $z_1,\ldots,
z_k\in \dC_+$ counting algebraic multiplicity $l_1, \ldots, l_k$
 are known as well as
$$b_1,a_1,b_2,\dots, a_{n-r-1}, b_{n-r},$$
 where
$$r=l_1+l_2+\cdots+
l_k.$$ Then $a_{n-r}, b_{n-r+1},\ldots,a_{n-1}, b_n$ are uniquely
determined.
\end{theorem}
\begin{proof} The sub-matrix $J_{[1,n-r]}$ is known. Therefore, the $m_-$-function $m_-(z,n-r+1)$ is known.
The unknown Jacobi sub-matrix $J_{[n-r+1, 1]}$ is self-adjoint and
 the corresponding Weyl function $m_+(z, n-r)$ belongs to the Herglotz-Nevanlinna class.
By Cramer's rule we have
that
\[
m_+(z,n-r)=\frac{L_{r-1}(z)}{L_r(z)},
\]
where $L_{r-1}(z)$ and $L_r(z)$ are polynomials with real
coefficients of degrees $r-1$ and $r$, respectively.
 Let
 $$\Phi(z):=a^2_{n-r}\,m_+(z,n-r).$$
By Theorem \ref{Green} we have 
\[
\Phi^{(p)}(z_j)=\left(\frac{1}{m_-(z,n-r+1)}\right)^{(p)}\Biggl
|_{z=z_j},\; j=1,\ldots,k, \; p=0,\ldots l_j.
\]
Therefore the values
\[
\Phi^{(p)}(z_j),\; j=1,\ldots,k, \; p=0,\ldots l_j
\]
are known.  Since $\Phi(\bar z)=\overline{\Phi(z)}$ and  $\IM
z_j>0$, the values
\[
\Phi^{(p)}(\overline{z_j})=\overline{\Phi^{(p)}(z_j)},\;
j=1,\ldots,k, \; p=0,\ldots l_j.
\]
are also known. Let $\wt\Phi(z)$ be a ratio of two polynomials with
real coefficients of degrees $r-1$ and $r$ such that
\[
\wt \Phi^{(p)}(z_j)=\Phi^{(p)}(z_j),\; \wt
\Phi^{(p)}(\bar{z_j})=\Phi^{(p)}(\bar{z_j}) j=1,\ldots,k, \;
p=0,\ldots l_j.
\]
Then for the function $\Psi(z)=\wt \Phi(z)-\Phi(z)$ we get
\[
\Psi^{(p)}(z_j)=\Psi^{(p)}(\overline{z_j})=0,\; j=1,\ldots,k, \;
p=0,\ldots l_j.
\]
It follows that
\[
\Psi(z)=\psi(z)\prod\limits_{j=1}^k(z-z_j)^{l_j}(z-\overline{z_j})^{l_j},
\]
where the function $\psi(z)$ has no zeroes at the points
$\{z_j\}_{j=1}^k$ and $\{\overline{z_j}\}_{j=1}^k$.
 The function $\Psi(z)$ is the ratio of two polynomials. The numerator is of degree less or equal
 $2r-1$.
 Since $r=l_1+\ldots+l_k$ we
get $\Psi(z)\equiv 0.$ So, $\Phi(z)$ is uniquely determined. Since
\[
a^2_{n-r}=-\lim\limits_{z\to\infty }z\Phi(z),
\]
we get the function $m_+(z,n-r)$ which uniquely determines the
self-adjoint Jacobi matrix $J_{[n-r+1,n]}$, i.e. the entries
$b_{n-r+1},\ldots,a_{n-1}, b_n$ are uniquely determined.
\end{proof}
\subsection{Extension and refinement of the Hochstadt and Gesztesy--Simon
uniqueness theorems for tri-diagonal matrices} \label{subsec55} Let
$J$ be the $n\times n$ tri-diagonal matrix of the form \eqref{1}
with the conditions
\begin{equation}
\label{trid} \left\{\begin{array}{l}\IM b_1>0, b_k=\bar b_k,\;
k=2,\ldots,n,\\
a_k > 0,\; k\ne p,\; a_p=0.
\end{array}\right.
\end{equation}
Recall (see Subsection \ref{subsec35}) that we consider $a$'s and
$b$'s as a single sequence $\{c_k\}_{k=1}^{2n-1}$, where
$c_{2k-1}=b_k$ and $c_{2k}=a_k$. The next theorem is the extension
and refinement of the Hochstadt theorem \cite{H3}.
\begin{theorem}
\label{Hh} Let $J$ be an $n\times n$ tri-diagonal matrix of the form
\eqref{1} with conditions \eqref{trid}. Suppose that $p$ non-real
eigenvalues $z_1,\ldots,z_p$ counting algebraic
multiplicity and distinct real eigenvalues \\
$z_{p+1},\ldots,z_n$
of $J$ are given as well as $c_{p+n+1},\ldots, c_{2n-1}.$ Then these
data uniquely determine entries $c_1,\ldots, c_{p+n}$.
\end{theorem}
\begin{proof}
Since $a_p=0$ the matrix $J$ takes the block form
\begin{equation}
\label{JJ} J=\begin{pmatrix}J_{11}&0\cr 0&J_{22} \end{pmatrix},
\end{equation}
where $J_{11}$ and $J_{22}$ are $p\times p$ and $(n-p)\times (n-p)$
tri-diagonal matrices of the form
\begin{equation}
\label{JJ1}
 J_{11}=\begin{pmatrix} b_1 & a_1 & 0    & 0 & \cdot &
\cdot & \cdot \\
a_1 & b_2 & a_2 & 0 & \cdot &
\cdot & \cdot \\
0    & a_2 & b_3 & a_3 & \cdot &
\cdot & \cdot \\
\cdot & \cdot & \cdot & \cdot &
\cdot & \cdot & \cdot \\
\cdot & \cdot & \cdot & \cdot &
\cdot & \cdot & a_{p-1} \\
\cdot & \cdot & \cdot & \cdot & 0 & a_{p-1} & b_p
\end{pmatrix},
\end{equation}
\begin{equation}
\label{JJ2}
 J_{22}=\begin{pmatrix} b_{p+1} & a_{p+1} & 0    & 0 &
\cdot &
\cdot & \cdot \\
a_{p+1} & b_{p+2} & a_{p+2} & 0 & \cdot &
\cdot & \cdot \\
\cdot & \cdot & \cdot & \cdot &
\cdot & \cdot & \cdot \\
\cdot & \cdot & \cdot & \cdot &
\cdot & \cdot & a_{n-1} \\
\cdot & \cdot & \cdot & \cdot & 0 & a_{n-1} & b_n
\end{pmatrix}.
\end{equation}
 The matrix $J_{22}$ is self-adjoint and therefore it has only real
distinct $n-p$ eigenvalues which are also eigenvalues of $J$. The
matrix $J_{11}$ is prime dissipative Jacobi matrix with a rank one
imaginary part. It follows that $J_{11}$ as well as $J$ has $p$
non-real eigenvalues counting multiplicity.
 By Theorem
\ref{TT1} the entries $c_1, c_2,\ldots, c_{2p-1}$ can be
reconstructed by means of the spectral data $z_1,\ldots z_p$. The
distinct real numbers $z_p,\ldots , z_n$ are eigenvalues of
self-adjoint Jacobi matrix $J_{22}$. By the Hochstadt theorem the
entries $c_{2p+1},\ldots, c_{2p+(n-p)}$ are determined uniquely by
$z_{p+1},\ldots,z_n$ and by the known entries $c_{p+n+1},\ldots,
c_{2n-1}.$
\end{proof}
Next is the extension and refinement of the Gesztesy--Simon theorem
\cite{GS}.
\begin{theorem}
\label{AGS} Let $J$ be an $n\times n$ tri-diagonal matrix of the
form \eqref{1} with conditions \eqref{trid}.
 Suppose that $p$ non-real eigenvalues
$z_1,\ldots,z_p$ counting algebraic multiplicity and distinct real
$j$ eigenvalues of $J$ are given as well as entries
\[
c_{2p+j+1},\ldots, c_{2n-1}.
\]
Then these data uniquely determine the matrix $J$.
\end{theorem}
\begin {proof}
From the proof of Theorem \ref{Hh} and conditions \eqref{trid} the
tri-diagonal matrix $J$ takes the form \eqref{JJ}. By Theorem
\ref{TT1} the entries $c_1, c_2,\ldots, c_{2p-1}$ can be
reconstructed uniquely by means of the non-real spectral data
$z_1,\ldots z_p.$ Using the Gesztesy--Simon uniqueness Theorem
\ref{GST} the $j$ distinct real eigenvalues and entries
$c_{2p+j+1},\ldots, c_{2n-1}$ uniquely determine the matrix
$J_{22}$. Therefore the given data uniquely determines matrix $J$.

\end{proof}
\begin{remark}
\label{RRRR}
 Consider more general tri-diagonal $n\times n$ dissipative matrix $J$ with
 $\ran J_I=\span\{\delta_1\}$
 of the form
 \[
J=\begin{pmatrix} b_1 & a_1 & 0    & 0 & \cdot &
\cdot & \cdot \\
\ovl{a}_1 & b_2 & a_2 & 0 & \cdot &
\cdot & \cdot \\
0    & \ovl{a}_2 & b_3 & a_3 & \cdot &
\cdot & \cdot \\
\cdot & \cdot & \cdot & \cdot &
\cdot & \cdot & \cdot \\
\cdot & \cdot & \cdot & \cdot &
\cdot & \cdot & a_{n-1} \\
\cdot & \cdot & \cdot & \cdot & 0 & \ovl{a}_{n-1} & b_n
\end{pmatrix}
 \]
 with the condition
\[
\IM b_1>0,\; b_k=\bar b_k,\; k=2,\ldots,n
\]
For such matrices the following statements are equivalent:
\begin{enumerate}
\item the matrix $J$ has $p$ non-real
eigenvalues counting their algebraic multiplicity and $p<n$;
\item $a_p=0$ and $a_k\ne 0$ for $k=1,\ldots, p-1.$
\end{enumerate}
\end{remark}
Actually, suppose that $a_p=0$ and all $a_k\ne 0$ for
$k=1,2,\ldots,p-1$. Then $J$ takes the  block form \eqref{JJ}, where
$J_{11}$ and $J_{22}$ are $p\times p$ and $(n-p)\times (n-p)$
tri-diagonal matrices given by \eqref{JJ1} and \eqref{JJ2}. The
matrix $J_{11}$ is prime dissipative Jacobi matrix with a rank one
imaginary part and $J_{22}$ is self-adjoint.  It follows that $J$
 has $p$ non-real eigenvalues counting algebraic
multiplicity.

Conversely, suppose that the matrix $J$ has $p$ non-real eigenvalues
counting algebraic multiplicity and $p<n$. Let $J$ be the
corresponding operator in $\dC^n$. Since $\ran
J_I=\span\{\delta_1\}$, the subspace
\[
\cH_1=\span\{\delta_1,J\delta_1,\ldots J^n\delta_1\}
\]
 reduces $J$, $J_{11}=J\uphar\cH_1$ is a prime dissipative
operator with a rank one imaginary part, $
J_{22}=J\uphar(\dC^n\ominus\cH_1)$ is a self-adjoint operator. Hence
$J_{11}$ has only non-real eigenvalues that coincide with non-real
eigenvalues of $J$. It follows that $\dim\cH_1=p$. From
\eqref{coord} it follows that $a_k\ne 0$ for $k=1,\ldots, p-1$, and
$a_p=0$.
\section{A non-self-adjoint analog of the Stone theorem}
\label{sec6}
  The next theorem is a non-self-adjoint analog of M.~
Stone's result (see Theorem \ref{Stone}).
\begin{theorem}
\label{Jacobimodel} Let $\cH$ be separable Hilbert space and let $A$
be a  bounded prime dissipative operator with a rank one imaginary
part acting in  $\cH$. Then there exists an orthonormal basis in
$\cH$ in which the matrix of the operator $A$ is a bounded
dissipative Jacobi matrix with conditions \eqref{333}.
\end{theorem}
\begin{proof} Let $\dim\cH=n$.
Because operator $A$ is prime, it has only $n$ non-real eigenvalues
$\{z_k\}\subset \dC_+$ counting their algebraic multiplicity. By
Theorem \ref{model} the characteristic function of $A$ takes the
form
\[
W(z)=\prod\limits_{k=1}^n\frac{z-{\bar z}_k}{z-z_k}.
\]
By Theorem \ref{TT1} there exists a Jacobi matrix
\[
\begin{pmatrix}
b_1 & a_1 & 0    & 0 & \cdot &
\cdot & \cdot \\
a_1 & b_2 & a_2 & 0 & \cdot &
\cdot & \cdot \\
0    & a_2 & b_3 & a_3 & \cdot &
\cdot & \cdot \\
\cdot & \cdot & \cdot & \cdot &
\cdot & \cdot & \cdot \\
\cdot & \cdot & \cdot & \cdot &
\cdot & \cdot & a_{n-1}\\
\cdot & \cdot & \cdot & \cdot & 0 & a_{n-1} & b_n
\end{pmatrix}
\]
with conditions \eqref{333} whose eigenvalues coincide with
$\{z_k\}$ and the corresponding characteristic function coincides
with $W(z)$. It follows that $A$ is unitarily equivalent to the
operator in $\dC^n$ determined by the matrix $J$.

Let $\dim\cH=\infty$ and let $g\in\cH$ such that $2A_Ih=(h,g)g,\;
h\in\cH.$ Define the characteristic function of $A$
\[
 W(z)=1-i\left((A-zI)^{-1}g,g\right), \; z\in\rho(A)
\]
 It follows that
\[
\lim\limits_{z\to\infty}z(W(z)-1)=i||g||^2.
\]
Let $c=||g||^2$. Define
\[
V(z)=i\,\frac{W(z)-1}{W(z)+1}=\frac{1}{2}\left((A_R-zI)^{-1}g,g\right)
\]
and let
\[
m(z)=\frac{2}{||g||^2}\,V(z).
\]
Then
\[
m(z)=\int\frac{d\rho(t)}{t-z},
\]
where the spectral measure has bounded support. Then there exists a
unique Jacobi matrix
\[
 H=\begin{pmatrix} b & a_1 & 0 &0   & 0 &
\cdot &
\cdot & \cdot \\
a_1 & b_2 & a_2 & 0 &0& \cdot &
\cdot & \cdot \\
0    & a_2 & b_3 & a_3 &0& \cdot &
\cdot & \cdot \\
\cdot & \cdot & \cdot & \cdot & \cdot & \cdot & \cdot &\cdot
\end{pmatrix}
\]
with real entries $b,b_2,b_3,\ldots$ and positive entries
$a_1,a_2,\ldots$ such that
$m(z)=\left((H-zI)^{-1}\delta_1,\delta_1\right).$  Note that
\eqref{Bound} holds because $H$ defines a bounded operator in
$l_2(\dN).$ Moreover, the entries of $H$ can be found by means of
continued fraction expansion
\[
m(z)= \frac{-1}{z-b}\;\raisebox{-3mm}{{\rm
+}}\;\frac{-a^2_1}{z-b_2}\;\raisebox{-3mm}{{\rm
+}}\;\frac{-a^2_2}{z-b_3} \;\raisebox{-3mm}{{\rm +}\,\ldots}\;
\raisebox{-3mm}{{\rm +}}\;\frac{-a^2_{n-1}}
{z-b_n}\;\raisebox{-3mm}{{\rm +}\,\dots}.
\]
 Let
\[
\begin{split}
&J=\begin{pmatrix} b+i||g||^2/2 & a_1 & 0 &0   & 0 & \cdot &
\cdot & \cdot \\
a_1 & b_2 & a_2 & 0 &0& \cdot &
\cdot & \cdot \\
0    & a_2 & b_3 & a_3 &0& \cdot &
\cdot & \cdot \\
\cdot & \cdot & \cdot & \cdot & \cdot & \cdot & \cdot &\cdot
\end{pmatrix}=\\
&\qquad\qquad=H+i\begin{pmatrix} ||g||^2/2 & 0 & 0    & 0 & \cdot &
\cdot & \cdot \\
0 & 0 & 0 & 0 & \cdot &
\cdot & \cdot \\
0    & 0 &0 & 0 & \cdot &
\cdot & \cdot \\
\cdot & \cdot & \cdot & \cdot &
\cdot & \cdot & \cdot \\
\end{pmatrix}.
\end{split}
\]
The characteristic function of $J$ is
\[
1-i ||g||^2\left((J-zI)^{-1}\delta_1,\delta_1\right).
\]
 Since $J_R=H$, and
$m(z)=\left((H-zI)^{-1}\delta_1,\delta_1\right),$ we get that the
characteristic function of $J$ coincides with $W(z)$. Because matrix
$J$ is
 prime, the operator $A$ is unitarily equivalent to $J$. Note that
the entries of $J$ can be also found using the continued fraction
expansion
\[
M(z)= \frac{-1}{z-b_1}\;\raisebox{-3mm}{{\rm
+}}\;\frac{-a^2_1}{z-b_2}\;\raisebox{-3mm}{{\rm
+}}\;\frac{-a^2_2}{z-b_3} \;\raisebox{-3mm}{{\rm +}\,\ldots}\;
\raisebox{-3mm}{{\rm +}}\;\frac{-a^2_{n-1}}
{z-b_n}\;\raisebox{-3mm}{{\rm +}\,\dots},
\]
where $M(z)=\frac{i}{\beta}(W(z)-1)$,
 and
 $\beta=\lim\limits_{z\to \infty}\left(iz(1-W(z))\right).$

\end{proof} From Theorem \ref{Jacobimodel} it follows that
the Jacobi matrices with conditions \eqref{Bound} and \eqref{333}
provide new models for bounded prime dissipative operators with rank
one imaginary part alternative to Livsic triangular models described
in Subsection \ref{TRMOD}.
\begin{remark}
From Theorems \ref{WeylFunction} and \ref{Jacobimodel} it follows
that the equivalent statements (iii), (iv), and (v) are also
equivalent to the statement: The function $M(z)$ is the Weyl
function of some bounded Jacobi matrix with conditions \eqref{333}.
\end{remark}
\begin{remark}
\label{GenJac} In the recent paper \cite{DD} it is proved that every
cyclic self-adjoint operator in a Pontryagin space is unitary
equivalent to some generalized Jacobi matrix \cite{KrLang}.
\end{remark}
The next theorem gives more information about the diagonal entries
of the model Jacobi matrix.
\begin{theorem}
\label{ZER} Let $A$ be a bounded prime dissipative operator with a
rank one imaginary part acting in a separable Hilbert space $\cH$.
Suppose that the characteristic function  $W(z)$ of $A$ possesses
the property
\begin{equation}
\label{WW}
 W(-z)=\frac{1}{W(z)},\; |z|>||A||.
\end{equation}
Then the corresponding to $A$ (by Theorem \ref{model}) bounded
Jacobi matrix of the form \eqref{22} with the conditions \eqref{333}
possesses the property
\[
\RE b_1=b_2=\cdots=0.
\]
\end{theorem}
\begin{proof} Since  (see \eqref{FLT})
\[
V(z)=i\,\frac{W(z)-1}{W(z)+1},\; z\in\rho(A)\cap\rho(A_R),
\]
we have  that in a neighborhood of infinity
\[
V(-z)=-V(z).
\]
Because
\[
m_{A_R}(z)=\frac{1}{l}\, V(z),
\]
where $l$ is the positive eigenvalue of $A_I$, we get that the Weyl
function of $A_R$ is odd. Since $m_{J_R}(z)=m_{A_R}(z)$, from
Proposition \ref{odd} it follows that $\RE b_1=b_2=\cdots=0$.
\end{proof}
\begin{corollary}
\label{CCCOL}  Let $A$ be a bounded prime dissipative operator with
a rank one imaginary part acting in a separable Hilbert space $\cH$.
Suppose that eigenvalues of $A$(counting algebraic multiplicity) are
symmetric with respect to the imaginary axis and $A$ has a complete
system of root subspaces. If $J$ is the corresponding to $A$ bounded
Jacobi matrix  with conditions \eqref{333}, then $J$ possesses the
property $\RE b_1=b_2=\cdots=0.$
\end{corollary}
\begin{proof} By Livsic  Theorem \ref{Livsic} the characteristic function
$W(z)$ takes the form
\[
W(z)=\prod\limits_{j=1}^N\left(\frac{z-\bar
z_j}{z-z_j}\right)^{l_j}\,\prod\limits_{j=1}^N\left(\frac{z+
z_j}{z+\bar
z_j}\right)^{l_j}\,\prod\limits_{k=1}^M\frac{z+i\eta_k}{z-i\eta_k},
\]
where $\{z_j\}_{j=1}^N$, $\{-\bar z_j\}_{j=1}^N$ $(\IM z_j>0,\;\RE
z_j< 0)$ are distinct eigenvalues of $A$, $\{l_j\}_{j=1}^N$ are
corresponding algebraic multiplicities, $\{i\eta_k\}_{k=1}^M$,
$\eta_k> 0$ are not necessarily distinct pure imaginary eigenvalues
of $A$, $N$ and $M$ are both finite when $\dim\cH<\infty$, and one
of them or both are infinite when $\dim\cH=\infty$, and
$M+2\,\sum\limits_{j=1}^N l_j=\dim\cH$. Note that if $M=0,$ (i.e.
there are no pure imaginary eigenvalues) then in the finite
dimensional $\cH$ its dimension is even. It is also possible that
$N=0$.

 It follows that
\[
W(-z)=\frac{1}{W(z)},\;|z|>||A||.
\]
By Theorem \ref{ZER} we obtain that the corresponding Jacobi matrix
possesses the property $\RE b_1=b_2=\cdots=0$.
\end{proof}
Now we want to mention some connections of Theorems \ref{TT1},
\ref{Jacobimodel} and Corollary \ref{CCCOL} with some result
established by K.~Veseli\'c in \cite{Ves1} and \cite{Ves2} related
to the spectral problems for the operator pencil of the form
\begin{equation}
\label{Pencil}
 z^2 M+z C+K,
\end{equation}
where $M$, $C$, and $K$ are self-adjoint $n\times n$ matrices with
real entries, $M$ and $K$ are positive definite, $C$ is positive
semi-definite. In \cite{Ves1} this problem is reduced to the
eigenvalue problem for some  $2n\times 2n$ matrix $A$ constructed by
means of $M,C$, and $K$. The inverse spectral problem for
\eqref{Pencil} is considered in \cite{Ves2} for the case
\begin{equation}
\label{MM}
 M={diag}(m_1,m_2,\ldots, m_n),
\end{equation}
\begin{equation}
\label{MK} K=\begin{pmatrix} k_1 & -k_1 & 0    & 0 & \cdot &
\cdot & \cdot \\
-k_1 &k_1+k_2 & -k_2 & 0 & \cdot &
\cdot & \cdot \\
\cdot & \cdot \\
\cdot & \cdot & \cdot & \cdot &
\cdot & \cdot & \cdot \\
\cdot & \cdot & \cdot & \cdot &
\cdot & \cdot & -k_{n-1} \\
\cdot & \cdot & \cdot & \cdot & 0 & -k_{n-1} & k_{n-1}+k_n
\end{pmatrix},
\end{equation}
and
\begin{equation}
\label{MC}
 C=\begin{pmatrix}
\gamma& 0 & 0 & 0 & \cdot &
\cdot & \cdot \\
0 & 0 & 0 & 0 & \cdot &
\cdot & \cdot \\
0    & 0 &0 & 0 & \cdot &
\cdot & \cdot \\
\cdot & \cdot & \cdot & \cdot &
\cdot & \cdot & \cdot \\
\cdot & \cdot & \cdot & \cdot &
\cdot & \cdot & \cdot \\
\cdot & \cdot & \cdot & \cdot & 0 & 0 & 0
\end{pmatrix},
\end{equation}
where $m_1=1,\;\{m_j\}, \{k_j\},$ and $\gamma$ are positive numbers.
In such situation the following result is established.
\begin{theorem}
\cite{Ves2}. \label{Ves}
 Prescribe $2n$ eigenvalues from the open
left half-plane which are symmetric (together with their algebraic
multiplicities) with respect to the real axis. Then there exist
unique matrices $M$, $K$ and $C$ of the form \eqref{MM}, \eqref{MK},
\eqref{MC}, respectively, such that
\begin{enumerate}
\def\labelenumi{\rm (\roman{enumi})}
\item $m_j>0,\; k_j>0,\; m_1=1,$
\item prescribed eigenvalues coincide with the spectrum of the operator pencil
\eqref{Pencil}.
\end{enumerate}
\end{theorem}
Note that in \cite{Ves1} and \cite{Ves2} it is supposed that the
total number of negative eigenvalues counting their algebraic
multiplicities is even.

As it was mentioned by K.~Veseli\'{c} (in private communication) the
matrix $A$ is unitarily similar to $2n\times 2n$ tri-diagonal matrix
$A'$ of the form
\[
A'=\begin{pmatrix} -\gamma & \alpha_1 & 0 & 0 & \cdot &
\cdot & \cdot \\
-\alpha_1 &0 & \beta_1 & 0 & \cdot &
\cdot & \cdot \\
0 & -\beta_1&0&\alpha_2&0&\cdot&\cdot \\
\cdot & \cdot & \cdot & \cdot &
\cdot & \cdot & \cdot \\
\cdot & \cdot & \cdot & \cdot &
\cdot & \cdot & \alpha_{n} \\
\cdot & \cdot & \cdot & \cdot & 0 & -\alpha_{n} &0
\end{pmatrix},
\]
with $\alpha_j=\sqrt{k_j/m_j},\; j=1,\ldots,n,\;
\beta_j=\sqrt{k_j/m_{j+1}},\;j=1,\ldots,n-1.$ The matrix $J'=-iA'$
is a dissipative tri-diagonal matrix with a rank one imaginary part.
Clearly, the matrix $J'$ is unitarily similar to $2n\times 2n$
Jacobi matrix
\[
J=\begin{pmatrix} i\gamma & \alpha_1 & 0 & 0 & \cdot &
\cdot & \cdot \\
\alpha_1 &0 & \beta_1 & 0 & \cdot &
\cdot & \cdot \\
0 & \beta_1&0&\alpha_2&0&\cdot&\cdot \\
\cdot & \cdot & \cdot & \cdot &
\cdot & \cdot & \cdot \\
\cdot & \cdot & \cdot & \cdot &
\cdot & \cdot & \alpha_{n} \\
\cdot & \cdot & \cdot & \cdot & 0 & \alpha_{n} &0\end{pmatrix}.
\]
It follows that the inverse spectral problem considered by
K.~Veseli\'c is reduced to the reconstruction problem of a
dissipative Jacobi matrix with a rank one imaginary part from its
eigenvalues. The symmetry property for the eigenvalues of the pencil
\eqref{Pencil} means that the eigenvalues of $J$ counting algebraic
multiplicity are symmetric with respect to the imaginary axis.
According to Theorem \ref{TT1} there exists a unique $2n\times 2n $
Jacobi matrix with entries satisfying conditions \eqref{333} whose
eigenvalues counting algebraic multiplicity are prescribed numbers
from the open upper half-plane and symmetric with respect to
imaginary axis. By Corollary \ref{CCCOL} the diagonal entries of the
reconstructed  Jacobi matrix possess the property
\[
\RE b_1=b_2=\ldots= b_{2n}=0.
\]
In conclusion we note that the proof of Theorem \ref{Ves} in
\cite{Ves2} differs from our proof of Theorem \ref{TT1} and does not
use the Livsic characteristic function.

\section{Dissipative Volterra operators with rank one imaginary parts and the corresponding semi-infinite
 Jacobi matrices}
 \label{sec7}
 Recall that a compact operator is called Volterra operator if its
 spectrum consists of one point (zero).
 \begin{theorem} \cite{L}
 \label{VolLiv}
Let $\cH$ be a separable Hilbert space. Let $A$ be a prime
dissipative Volterra operator with a rank one imaginary part and let
 \begin{equation}
\label{IMM}
  A_I h=l(h,e)e
\end{equation}
for every $h$, where $l>0$, $||e||=1$. Then $A$ is unitarily
equivalent to the integration operator
\[
\c(Ff)(x)=2i\int\limits_{x}^lf(t)dt
\]
in the Hilbert space $\cL_2[0,l]$.
\end{theorem}
In connection with Theorem \ref{Jacobimodel} we will find a Jacobi
matrix corresponding to Volterra operator.
\begin{theorem}
\label{CANON} Let $A$ be a prime dissipative Volterra operator with
a rank one imaginary part and let $A_I h=l(h,e)e$ for every $h$,
where $l>0$, $||e||=1$. Then $A$ is unitarily equivalent to the
operator in $l_2(\dN)$ determined by the dissipative Jacobi matrix
\begin{equation}
\label{CAN} J=\begin{pmatrix} il & \frac{l}{\sqrt{1\cdot 3}} & 0 &0
& 0 &
\cdot & \cdot &\cdot&\cdot&\cdot\\
\frac{l}{\sqrt{1\cdot 3}}& 0 & \frac{l}{\sqrt{3\cdot 5}} & 0 &0&
\cdot &
\cdot & \cdot&\cdot \\
0    & \frac{l}{\sqrt{3\cdot 5}} & 0 & \frac{l}{\sqrt{5\cdot 7}} &0&
\cdot & \cdot &\cdot&\cdot&\cdot\\
\cdot & \cdot & \cdot & \cdot & \cdot & \cdot &\cdot\\
 \cdot & \cdot &
\cdot & \cdot & \cdot & \cdot &\frac{l}{\sqrt{(2n-3)
(2n-1)}}&0& \frac{l}{\sqrt{2n-1) (2n+1)}}&\cdot\\
\cdot & \cdot & \cdot & \cdot & \cdot & \cdot &\cdot&\cdot&\cdot&\cdot\\
\end{pmatrix}.
\end{equation}
\end{theorem}
\begin{proof}
As is well known \cite{L}, \cite{Br} that if \eqref{IMM} holds, then
the characteristic function $W(z)=1-2il\left((A-zI)^{-1}e,e\right)$
of $A$ takes the form
\[
W(z)=\exp\left(\frac{2il}{z}\right).
\]
It follows that
\[
V(z)={l}\left((A_R-zI)^{-1}e,e\right)=i\,\frac{W(z)-1}{W(z)+1}=-\tan\left(\frac{l}{z}\right).
\]
By Proposition \ref{SS} the vector $e$ is a cyclic for $A_R$,
therefore $m(z)=\left((A_R-zI)^{-1}e,e\right)$ is the Weyl function
of $A_R$. We get
\[
m(z)=-\frac{1}{l}\tan\left(\frac{l}{z}\right).
\]
 We consider
for simplicity the case $l=1$. In this case
$m(z)=-\tan\left(\frac{1}{z}\right).$ We have to find a self-adjoint
semi-infinite Jacobi matrix whose the Weyl  function is
$-\tan\left(\frac{1}{z}\right).$ This function has the following
continued fraction expansion \cite{KK}
\[
\tan
x=\frac{x}{1}\;\raisebox{-3mm}{--}\;\frac{x^2}{3}\;\raisebox{-3mm}{--}\;\frac{x^2}{5}
\;\raisebox{-3mm}{--}\;\frac{x^2}{7}\;\raisebox{-3mm}{--}\;\frac{x^2}{9}\;\;
\raisebox{-3mm}{--\,\ldots}\; \raisebox{-3mm}{--}\;\frac{x^2}
{2n+1}\;\raisebox{-3mm}{--\,\dots}.
\]
 Comparing with \eqref{CF} and taking into account that $x=1/z$ we get
\[
b_1=b_2=\ldots=0,\; a^2_n=\frac{1}{(2n-1)(2n+1)},\; n=1,2,\ldots.
\]
It follows that the Jacobi matrix with the Weyl function
$-\tan\left(\frac{1}{z}\right)$ is the following
\[
H_0=\begin{pmatrix} 0 & \frac{1}{\sqrt{1\cdot 3}} & 0 &0 & 0 &
\cdot & \cdot &\cdot&\cdot&\cdot\\
\frac{1}{\sqrt{1\cdot 3}}& 0 & \frac{1}{\sqrt{3\cdot 5}} & 0 &0&
\cdot &
\cdot & \cdot&\cdot \\
0    & \frac{1}{\sqrt{3\cdot 5}} & 0 & \frac{1}{\sqrt{5\cdot 7}} &0&
\cdot & \cdot &\cdot&\cdot&\cdot\\
\cdot & \cdot & \cdot & \cdot & \cdot & \cdot &\cdot\\
 \cdot & \cdot &
\cdot & \cdot & \cdot & \cdot &\frac{1}{\sqrt{(2n-3)
(2n-1)}}&0& \frac{1}{\sqrt{2n-1) (2n+1)}}&\cdot\\
\cdot & \cdot & \cdot & \cdot & \cdot & \cdot &\cdot&\cdot&\cdot&\cdot\\
\end{pmatrix}.
\]
Let
\[
J_0=\begin{pmatrix} i & \frac{1}{\sqrt{1\cdot 3}} & 0 &0 & 0 &
\cdot & \cdot &\cdot&\cdot&\cdot\\
\frac{1}{\sqrt{1\cdot 3}}& 0 & \frac{1}{\sqrt{3\cdot 5}} & 0 &0&
\cdot &
\cdot & \cdot&\cdot \\
0    & \frac{1}{\sqrt{3\cdot 5}} & 0 & \frac{1}{\sqrt{5\cdot 7}} &0&
\cdot & \cdot &\cdot&\cdot&\cdot\\
\cdot & \cdot & \cdot & \cdot & \cdot & \cdot &\cdot\\
\cdot & \cdot & \cdot & \cdot & \cdot & \cdot &\cdot&\cdot&\cdot&\cdot\\
 \cdot & \cdot &
\cdot & \cdot & \cdot & \cdot &\frac{1}{\sqrt{(2n-3)
(2n-1)}}&0& \frac{1}{\sqrt{2n-1) (2n+1)}}&\cdot\\
\cdot & \cdot & \cdot & \cdot & \cdot & \cdot
&\cdot&\cdot&\cdot&\cdot
\end{pmatrix}.
\]
Then the characteristic function of $J_0$ is
$\exp\left(\frac{2i}{z}\right)$. So the operator $J$ is unitary
equivalent to $A$ with $l=1.$ For arbitrary $l$ the Jacobi matrix
$J=lJ_0$ takes the form \eqref{CAN} and its characteristic function
is $\exp\left(\frac{2il}{z}\right)$.
\end{proof}
\begin{remark}
\label{REM1} The real part $\cF_R$ of the integration operator $\c(F
f)(x)=2i\int\limits_{x}^lf(t)dt$ takes the form
\[
\c(F_Rf)(x)=i\left(\int\limits_{x}^lf(t)dt-\int\limits_{0}^xf(t)dt\right).
\]
Its eigenvalues are
\[
\frac{2l}{(2k+1)\pi},\; k\in\dZ.
\]
From Theorem \ref{CANON} it follows that these numbers are
eigenvalues of the matrix
\begin{equation}
\label{H}
 H=\begin{pmatrix} 0 & \frac{l}{\sqrt{1\cdot 3}} & 0 &0 & 0
&
\cdot & \cdot &\cdot&\cdot&\cdot\\
\frac{l}{\sqrt{1\cdot 3}}& 0 & \frac{l}{\sqrt{3\cdot 5}} & 0 &0&
\cdot &
\cdot & \cdot&\cdot \\
0    & \frac{l}{\sqrt{3\cdot 5}} & 0 & \frac{l}{\sqrt{5\cdot 7}} &0&
\cdot & \cdot &\cdot&\cdot&\cdot\\
\cdot & \cdot & \cdot & \cdot & \cdot & \cdot &\cdot\\
 \cdot & \cdot &
\cdot & \cdot & \cdot & \cdot &\frac{l}{\sqrt{(2n-3)
(2n-1)}}&0& \frac{l}{\sqrt{2n-1) (2n+1)}}&\cdot\\
\cdot & \cdot & \cdot & \cdot & \cdot & \cdot &\cdot&\cdot&\cdot&\cdot\\
\end{pmatrix}.
\end{equation}
\end{remark}
\begin{remark}
\label{REM2}
 It is known \cite{Akh}, \cite{S} that there are
relations between the moments $\{\gamma_k=\int t^k d\rho(t)\}$,
$\gamma_0=1$ of a measure $d\rho(t)$ and the entries $\{b_k\}$,
$\{a_k\}$ of the corresponding self-adjoint Jacobi matrix  \[
J=\begin{pmatrix} b_1 & a_1 & 0 &0   & 0 & \cdot &
\cdot & \cdot \\
a_1 & b_2 & a_2 & 0 &0& \cdot &
\cdot & \cdot \\
0    & a_2 & b_3 & a_3 &0& \cdot &
\cdot & \cdot \\
\cdot & \cdot & \cdot & \cdot & \cdot & \cdot & \cdot &\cdot
\end{pmatrix}.
\]
More precisely, let
\[
\begin{split}
& h_0=1,\; h_k=\det\begin{pmatrix} \gamma_0 &
\gamma_1&\cdot&\cdot&\cdot &\gamma_{k-1}
 \\
\gamma_1 &\gamma_2&\cdot&\cdot&\cdot&\gamma_{k}
 \\
\cdot & \cdot & \cdot&\cdot&\cdot&\cdot  \\
\gamma_{k-1}&\gamma_k &\cdot&\cdot&\cdot&\gamma_{2k-2}
\end{pmatrix},\; k\ge 1,\\
&\wt h_k= \det\begin{pmatrix} \gamma_0 & \gamma_1&\cdot&\cdot&\cdot
&\gamma_{k-2}&\gamma_k
 \\
\gamma_1 &\gamma_2&\cdot&\cdot&\cdot&\gamma_{k-1}&\gamma_{k+1}
 \\
\cdot & \cdot & \cdot&\cdot&\cdot&\cdot& \cdot  \\
\gamma_{k-1}&\gamma_k
&\cdot&\cdot&\cdot&\gamma_{2k-3}&\gamma_{2k-1},
\end{pmatrix}\; k\ge 2.
\end{split}
\]
Then
 \[
a_k=\left(\frac {h_{k-1}\,
h_{k+1}}{h^2_k}\right)^{1/2},\;\;b_1=\gamma_1,\; \sum\limits_{j=1}^k
b_j=\frac{\wt h_{k}}{h_k}.
\]
Since
\[
\left((H-zI)^{-1}\delta_1,\delta_1\right)=-\frac{1}{l}\tan\left(\frac{l}{z}\right),
\]
where $H$ is given by \eqref{H}, the function
$\tan\left(\frac{1}{z}\right)$ has the following Taylor expansion at
infinity
\[
\tan\left(\frac{1}{z}\right)=\sum\limits_{n=1}^\infty\frac{2^{2n}(2^{2n}-1)}{(2n)!}\,\frac{B_{2n}}{z^{2n-1}},
\]
where $\{B_{2n}\}$ are Bernoulli numbers \cite{KK}, $B_2=1/6,\;
B_{4}=-1/30,\; B_{6}=1/42,\ldots$, we get that
\[
\gamma_{2k}=-\frac{l^{2(k-1)}2^{2k}(2^{2k}-1)B_{2k}}{(2k)!},\;\,
\gamma_{2k-1}=0,\; k=1,2,\ldots.
\]
The formula for $\{a_k\}$  allows to get some relations between
Bernoulli numbers.
\end{remark}
\begin{theorem}
\label{CCCC} Let $J_t$ be a Jacobi matrix of the form
\[
J_t=\begin{pmatrix} t & \frac{l}{\sqrt{1\cdot 3}} & 0 &0 & 0 &
\cdot & \cdot &\cdot&\cdot&\cdot\\
\frac{l}{\sqrt{1\cdot 3}}& 0 & \frac{l}{\sqrt{3\cdot 5}} & 0 &0&
\cdot &
\cdot & \cdot&\cdot \\
0    & \frac{l}{\sqrt{3\cdot 5}} & 0 & \frac{l}{\sqrt{5\cdot 7}} &0&
\cdot & \cdot &\cdot&\cdot&\cdot\\
\cdot & \cdot & \cdot & \cdot & \cdot & \cdot &\cdot\\
 \cdot & \cdot &
\cdot & \cdot & \cdot & \cdot &\frac{l}{\sqrt{(2n-3)
(2n-1)}}&0& \frac{l}{\sqrt{2n-1) (2n+1)}}&\cdot\\
\cdot & \cdot & \cdot & \cdot & \cdot & \cdot &\cdot&\cdot&\cdot&\cdot\\
\end{pmatrix},
\]
where $t$ is a complex number and $l$ is a positive number. Then the
corresponding Jacobi operator for all $t$ but two $t=\pm i\,l$ has a
complete system of eigenvectors.
\end{theorem}
\begin{proof}
The matrix $J_t$ has the representation
\[
J_t=H+t(\cdot,\delta_1)\delta_1,
\]
where $H$ is given by \eqref{H}. Suppose that $\IM t\ne 0$. Then the
non-real eigenvalues of $J_t$ are solutions of the equation
\[
1+t\left((H-zI)^{-1}\delta_1,\delta_1\right)=0.
\]
Because
\[
\left((H-zI)^{-1}\delta_1,\delta_1\right)=-\frac{1}{l}\tan\left(\frac{l}{z}\right),
\]
we get the equation
\[
\tan\left(\frac{l}{z}\right)=\frac{l}{t}.
\]
This equation has the following solutions:
\[
z_k=\frac{2l\left(\arg\left(\displaystyle\frac
{t+il}{t-il}\right)+2\pi
k+i\,\ln\left|\displaystyle\frac{t+il}{t-il}\right|\right)}
{\left(\arg\left(\displaystyle\frac{t+il}{t-il}\right)+2\pi
k\right)^2+ \ln^2\left|\displaystyle\frac{t+il}{t-il}\right|},\;
k\in\dZ,\; t\ne\pm il.
\]
 Let
\[
x=\ln\left|\frac{t+il}{t-il}\right|,\;
y=\arg\left(\frac{t+il}{t-il}\right),\; t\ne\pm il.
\]
Then
\[
\sum\limits_{k=-\infty}^{\infty}\IM
z_k=l\,\sum\limits_{k=-\infty}^\infty\frac{2x}{x^2+(y+2\pi k)^2}\,.
\]
It is well known that
\[
\sum\limits_{k=-\infty}^\infty\frac{2\,x}{x^2+(y+2\pi
k)^2}=-\frac{\sinh x}{\cos y-\cosh x}\,.
\]
We have
\[
\begin{split}
&\sinh x=\frac{|t+il|^2-|t-il|^2}{2|t-il||t+il|},\\
&\cos
y=\cos\left(\arg\frac{t+il}{t-il}\right)=\frac{\RE\displaystyle
\frac{t+il}{t-il}}{\left|\displaystyle\frac{t+il}{t-il}\right|}=\RE\frac{(t+il)(\bar
t+il)}{|t-il||t+il|}\,,\\
&\cosh x=\frac{|t+il|^2+|t-il|^2}{2|t-il||t+il|}\,,\\
&-\frac{\sinh x}{\cos y-\cosh
x}=-\frac{|t+il|^2-|t-il|^2}{2\RE\left((t+il)(\bar
t+il)\right)-|t+il|^2-|t-il|^2}=\\
&=\frac{|t+il|^2-|t-il|^2}{4l^2}=\frac{\IM t}{l}.
\end{split}
\]
Hence
\[
\sum\limits_{k=-\infty}^{\infty}\IM z_k=\IM t.
\]
Since in this case the eigenvalues of $J_t$ have the algebraic
multiplicities one, according  to the Livsic Theorem \ref{Livsic1}
the operator $J_t$  has complete system of eigenvectors.

If $\IM t=0$, then the operator $J_t$ is self-adjoint and compact, its
eigenvalues are real numbers
\[
z_k=\frac{2l} {\arg\left(\displaystyle\frac{t+il}{t-il}\right)+2\pi
k}\,,\; k\in\dZ,
\]
and the corresponding eigenvectors form a complete system.
\end{proof}
{\bf Acknowledgement} We are thankful to Fritz Gesztesy, Konstantin
Makarov and Kresimir Veseli\'c for valuable discussions.


\begin{thebibliography}{33}

\bibitem{Akh}
N. I.~Akhiezer, The Classical Moment Problem, Oliver and Boyd,
Edinburgh, 1965.

\bibitem{AkhGl}
N. I.~Akhiezer and I. M.~Glazman, Theory of Linear Operators in
Hilbert spaces, volume I, II, Pitman, Boston 1981.


\bibitem{ArHdeS}
Yu. M.~Arlinski\u{i}, S.~Hassi, and H. S. V.~de Snoo,
\textit{$Q$-functions of quasi-selfadjoint contractions}, Operator
Theory: Advances and Applications, 163 (2006), 23--54.

\bibitem{ArTs2}
Yu. M.~Arlinski\u{\i} and E. R.~Tsekanovski\u{\i},
\textit{Quasi-selfadjoint contractive extensions of Hermitian
contractions}, Teor. Funkts., Funkts. Anal. Prilozhen, 50 (1988),
9--16.

\bibitem{Atk}
F.~Atkinson, Discrete an Continuous Boundary Problems, Academic
Press, New-York, 1964.

\bibitem{Bec1}
B.~Beckerman, \textit{Complex Jacobi matrices}, J. Comput. Appl.
Math. 127 (2001), No. 1-2, 17--65.



\bibitem{Bec}
B.~Beckerman, \textit{On the classificqtion of the spectrum of
second order Difference operator}, Math. Nachr., 916 (2000), 45--59.


\bibitem{Ber}
Yu. M.~Berezanski\u{i}, Expansions in Eighenfunctions of
Self-Adjoint Operators, Transl. Math. Mono. 17, Amer. Math. Soc.,
Providence, RI, 1968.

\bibitem{Ber1}
Yu. M.~Berezanski\u{i}, \textit {Integration of nonlinear difference
equation by the inverse spectral method}, Soviet Math. Doklady 31
(1985), No. 2, 264--267.

\bibitem{Br}
M. S.~Brodski\u{\i}, Triangular and Jordan Representations of Linear
Operators, Nauka, Moscow, 1969. (Russian).

\bibitem{BrL}
M. S.~Brodski\u{\i} and M. S.~Livsic, \textit{Spectral analysis of
nonself-adjoint operators and intermediate systems}, Uspekhi Mat.
Nauk, 13 (1958), No.1, 3--85. (Russian).

\bibitem{Ch}
S.~Chung-Tsun, \textit {Some inverse problems on Jacobi matrices},
 Inverse Problems 20 (2004), no. 2, 589-600.



\bibitem{DN}
P.~Deift and T.~Nanda, \textit {On the determination of a
tridiagonal matrix from its spectrum and submatrix}, Linear Algebra
and Appl. 60 (1984), 43--55.

\bibitem{DD}
M.~Derevyagin and V.~Derkach, \textit {Spectral problems for
generalized Jacobi matrices}, Linear Algebra and Appl. 382 (2004),
1--24.


\bibitem{GS}
F.~Gesztesy and B.~Simon, \textit{M-functions and inverse spectral
analysis for finite and semifinite Jacobi matrices}, J. Anal. Math.
73 (1997), 267--297.

\bibitem{GS1}
F.~Gesztesy and B.~Simon, \textit{Uniqueness theorems in inverse
spectral theory for one-dimensional Schr\"odinger operators}, Trans.
Amer. Math. Soc., 348 (1996), 349--373.

\bibitem{GL}
I. M.~Glazman and Yu. I.~Lyubich, Finite-Dimensional Linear
Analysis, Moscow, Nauka, 1969. (Russian) Dear Kolya, Dear Pavel,

\bibitem{G}
P.~Gibson, \textit {Inverse spectral theory of finite Jacobi
matrices}, Trans.Amer.Math.Soc., 354 (2002), no. 12, 4703--4749.

\bibitem{H1}
H.~Hochstadt, \textit{On some inverse problems in matrix theory},
Archiv der Math., 18 (1967), 201--207.

\bibitem{H2}
H.~Hochstadt, \textit{On the construction of a Jacobi matrix from
 spectral data}, Lin. Algebra and Appl., 8 (1974), 435--446.


\bibitem{H3}
H.~Hochstadt, \textit{On the construction of a Jacobi matrix from
 mixed given data}, Lin. Algebra and Appl., 28 (1979), 113--115.

\bibitem{HL}
H.~Hochstadt and B.~Lieberman, \textit{An inverse Sturm--Lioville
problem with mixed given data}, SIAM J. Appl. Math., 34 (1978),
676--680.

\bibitem{JTh}
W. B.~Jones and W. J.~Thron, Continued Fractions. Analytic Theory
and Applications, Addison-Wesley Publishing Company (1980).


\bibitem{KK}
 A.~Korn and M.~Korn, Mathematical Hanbook for Scientists and
Engineers, second edition, McGraw-Hill, 1968.

\bibitem{Kr}
M. G.~Kre\u{\i}n, \textit{The theory of selfadjoint extensions of
semibounded Hermitian operators and its applications}, Mat. Sb., 20
(1947) 431--495; 21 (1947), 365--404.

\bibitem {KrLang}
M.G.~Kre\u{\i}n and H.~Langer, \textit {On some extension problems
which are closely connected with the theory of hermitian operators
in a space $\Pi_\kappa$ III. Indefinite analogues of the Hamburger
and Stieltjes moment problems}, Part I, Beitr. Anal. 14 (1979),
25--40.

\bibitem{L}
M. S.~Livsic, \textit{On a spectral decomposition of linear
nonself-adjoint operator}, Amer. Math. Soc. Transl. (2) 5, 1957,
67--114.

\bibitem{LYan}
M. S.~Livsic and A. A.~Yantsevich, Operator Colligations in Hilbert
Spaces, Kharkov University, Kharkov, 1971. (Russian)

\bibitem{S}
B.~Simon, \textit{The classical moment problem as a self-adjoint
finite difference operator}, Advances in Math., 137 (1998), 82--203.

\bibitem{Ves1}
K.~Veseli\'{c}, \textit{On linear vibrational systems with one
dimensional damping}, Appl.Anal., 29 (1988), 1--18.

\bibitem{Ves2}
K.~Veseli\'{c}, \textit{On linear vibrational systems with one
dimensional damping II}, Integral Eq. Oper. Theory, 13 (1990),
883--897.


\bibitem{Stone}
M.~Stone, Linear Transformations in Hilbert Space and Their
Applications to Analysis, New-York, 1932.

\bibitem{Tes}
G.~Teschl, Jacobi Operators and Completely Integrable Nonlinear
Lattices, Math. Surv. and Mon. 72, Amer. Math. Soc., Rhode Island,
2000.

\bibitem{Wall}
H. S.~Wall, Analytic Theory of Continued Fractions, Chelsea, Bronx
NY (1973).
\end{thebibliography}
\end{document}